\documentclass[11pt]{article}

\usepackage[T1]{fontenc}
\usepackage[french]{babel}

\usepackage{amsmath,latexsym,enumerate,amsthm,amssymb}

\newcommand{\al}{\alpha}
\newcommand{\be}{\beta}
\newcommand{\ga}{\gamma}
\newcommand{\ep}{\epsilon}
\newcommand{\la}{\lambda}
\newcommand{\om}{\omega}
\newcommand{\si}{\sigma}
\renewcommand{\th}{\theta}

\newcommand{\C}{\mathbb{C}}

\newcommand{\N}{\mathbb{N}}
\newcommand{\Q}{\mathbb{Q}}
\newcommand{\R}{\mathbb{R}}
\newcommand{\T}{\mathbb{T}}
\newcommand{\Z}{\mathbb{Z}}

\newcommand{\zerob}{\boldsymbol{0}}
\newcommand{\unb}{\boldsymbol{1}}

\newcommand{\eb}{\mathbf{e}}
\newcommand{\jb}{\mathbf{j}}
\newcommand{\kb}{\mathbf{k}}
\newcommand{\lb}{\boldsymbol{\ell}}
\newcommand{\mb}{\mathbf{m}}
\newcommand{\rb}{\mathbf{r}}
\renewcommand{\sb}{\mathbf{s}}
\newcommand{\ub}{\mathbf{u}}
\newcommand{\vb}{\mathbf{v}}
\newcommand{\xb}{\mathbf{x}}
\newcommand{\yb}{\mathbf{y}}
\newcommand{\zb}{\mathbf{z}}

\newcommand{\Xb}{\mathbf{X}}

\newcommand{\ab}{\boldsymbol{\alpha}}
\newcommand{\bb}{\boldsymbol{\beta}}
\newcommand{\gb}{\boldsymbol{\gamma}}

\newcommand{\mub}{\boldsymbol{\mu}}
\newcommand{\sib}{\boldsymbol{\sigma}}
\newcommand{\taub}{\boldsymbol{\tau}}

\newcommand{\p}{\partial}

\newcommand{\ds}{\displaystyle}

\newcommand{\br}{{\cal B}}
\newcommand{\Cr}{{\cal C}}
\newcommand{\er}{{\cal E}}

\newcounter{nums}[section]

\newtheorem{cor}[nums]{Corollaire}
\newtheorem{déf}[nums]{Définition}
\newtheorem{ex}[nums]{Exemple}
\newtheorem{lemme}[nums]{Lemme}
\newtheorem{notation}[nums]{Notation}
\newtheorem{prop}[nums]{Proposition}
\newtheorem{rem}[nums]{Remarque}
\newtheorem{theo}{Théorème}

\begin{document}

\begin{center}
Valeurs aux T-uplets d'entiers négatifs \\
 de séries z\^etas multivariables \\ associées \\
à des polynômes de plusieurs variables.  \\ 
\ \\
Marc de Crisenoy \\
adresse électronique: mdecrise@math.unicaen.fr \\
\ \\
Laboratoire de Mathématiques Nicolas Oresme \\
CNRS UMR 6139 \\
LMNO \\
Université de Caen BP 5186 \\
F 14032 Caen Cedex \\

\end{center}

\part*{Introduction}

   Soient $ Q,P_1,...,P_T \in \R[X_1,...,X_N] $ et $ \mu_1,...,\mu_N $ des nombres complexes de module $ 1 $. \\
On considère la série de Dirichlet généralisée suivante: 
$$ Z(Q,P_1,...,P_T,\mu_1,...,\mu_N,s_1,...,s_T) = \sum_{m_1 \geq 1,...,m_N \geq 1} \frac{(\prod_{n=1}^N \mu_n^{m_n}) Q(m_1,...,m_N)}{\prod_{t=1}^T P_t(m_1,...,m_N)^{s_t}} $$
o\`u $ (s_1,...,s_T) \in \C^T $. \\
Le prolongement de ces séries a été étudié successivement par Mahler (\cite{mahler}), Mellin (\cite{mellin}), Cassou-Noguès (\cite{cassou 82}), Sargos (\cite{sargos 84}), Lichtin (\cite{lichtin 88}) et Essouabri (\cite{essoua}).

  Une simple adaptation du résultat de \cite{essoua} permet de voir que si $ P_1,...,P_T $ vérifient l'hypothèse probablement optimale $ H_0S $ (voir ci-dessous), alors cette série se prolonge méromorphiquement à $ \C^T $. On s'attend à ce que, lorsque la série est réellement tordue (c'est à dire lorsque $ \mu_1,...,\mu_N $ sont tous différents de $ 1 $) le prolongement soit holomorphe. Comme nous le montrerons sur un exemple, ce n'est pas toujours le cas. Dans ce travail nous introduisons une classe de polyn\^omes  contenant strictement celle des polyn\^omes à coefficients positifs et contenu dans celle des polyn\^omes vérifiant $ H_0S $. Nous montrons que dans cette classe le prolongement de $ Z $ est {\bf holomorphe} sur $ \C^T $. L'utilisation de ces séries multivariables (ie $ T $ quelconque) fournit le cadre naturel d'un {\bf lemme d'échange} crucial. De ce lemme d'échange on déduit le principal résultat de ce travail à savoir des formules simples et explicites pour les valeurs aux points $ \sb = (-k_1,...,-k_T) \in (-\N)^T $ ($ T $-uplet d'entiers négatifs). Nous transformons alors ces formules pour retrouver celles de Cassou-Noguès permettant de réaliser l'interpolation p-adique. Rappelons qu'elle en avait déduit l'existence des fonctions zêtas $p-$adique associées aux corps de nombres totalement réels. Ce résultat avait aussi été obtenu indépendamment par Barsky (\cite{barsky}) et Deligne-Ribet (\cite{deligne; ribet}). \\
L'étude des valeurs aux entiers négatifs de séries de Dirichlet associées à des polyn\^omes de plusieurs variables est un ancien problème lié en particulier aux propriétés de divers fonctions zéta intervenant dans l'arithmétique des corps de nombres (fonction zéta de Dedekind...cf Shintani (\cite{shintani})). Les résultats les plus importants en lien avec notre travail sont ceux obtenus par Pierrette Cassou-Noguès(\cite{cassou 79} et \cite{cassou 82}) dans le cas d'un polyn\^ome à coefficients positifs ( $ T = 1 $ ). Citons aussi le travail de Kwang-Wu Chen et Minking Eie (\cite{chen; eie}) qui ont obtenu sous les m\^emes hypothèses que Pierrette Cassou-Noguès et par des méthodes semblables aux siennes des formules très simples pour les valeurs aux entiers négatifs $ -k $. \\
Dans notre travail nous obtenons des formules aussi simples que celles de Kwang-Wu Chen et Minking Eie mais pour une classe plus générale de polyn\^omes (classe HDF) et pour des séries associées à plusieurs polyn\^omes ( $ T $ quelconque). Par ailleurs nos méthodes sont radicalement différentes de celles utilisées par Cassou-Noguès et Chen-Eie; nous espérons que la méthode du lemme d'échange donne une meilleure compréhension de la nature de ces formules. \\
 Signalons enfin que lorsque les séries ne sont pas (forcément) tordues le cas très particulier des formes linéaires a été étudié par divers auteurs (Akiyama, Egami et Tanigawa dans \cite{akiyama; egami; tanigawa}; Akiyama et Ishikawa dans \cite{akiyama; ishikawa}; Akiyama et Tanigawa dans \cite{akiyama; tanigawa}; Arakawa et Kaneko dans \cite{arakawa; kaneko}; Egami et Matsumoto dans \cite{egami; matsumoto}; Zhao dans \cite{zhao} ...).

\part*{Plan}

\noindent Enoncés des principaux résultats. \\
Quelques remarques. \\
Des lemmes sur les polyn\^omes à plusieurs variables. \\
Prolongement holomorphe des intégrales $ Y $. \\
Domaine de convergence de $ Z $. \\
Représentation intégrale. \\
Preuve du théorème \ref{existence du prolongement holomorphe de Z}. \\
Lemme d'échange et valeurs de $ Z $ aux points de $ (-\N)^T $. \\
Une formule pour les valeurs de $ \zeta_{\mu} $ aux entiers négatifs. \\
Preuve des théorèmes 4 et 5. \\

\part*{Enoncés des principaux résultats}

\noindent {\bf Convention}: dans tout ce travail on dira d'une série à $ N \geq 1 $ variables qu'elle est convergente lorsqu'elle est sommable (au sens des familles sommables). \\
En particulier, une série à une variable est dite convergente lorsqu'elle est absolument convergente. \\

\begin{notation}
On note $ \T = \{ z \in \C \ | \ |z| = 1 \} $ et $ J = [1, + \infty[ $. \\
\end{notation}

\begin{déf}
Soit $ \mu \in \T $. \\
On pose $ \ds \zeta_{\mu}(s) = \sum_{m \geq 1} \frac{\mu^m}{m^s} $. \\
\end{déf}

\begin{déf}
Soient $ Q,P_1,...,P_T \in \R[X_1,...,X_N] $ tels que $ \forall 1 \leq t \leq T \ \ \forall \xb \in J^N \ P_t(\xb) > 0 $. \\
Soit $ \mub \in {\T}^N $. \\
On pose :
$$ Z(Q,P_1,...,P_T,\mub,s_1,...,s_T) = \sum_{\mb \in \N^{*N}} \mub^{\mb} Q(\mb) \prod_{t=1}^T P_t(\mb)^{- s_t} $$
\end{déf}

\begin{rem}
Pour $ \mu \in \T $ on a: $ Z(1,X,\mu,\cdot) = \zeta_{\mu} $. \\
\end{rem}

Introduisons une nouvelle classe de polyn\^omes: \\

\begin{déf}
soit $ P \in \R[X_1,...,X_N] $ \\ 
On dit que $ P $ vérifie l'hypothèse raisonnable (abrégée en HDF dans toute la suite) si: \\
$ \forall \xb \in J^N \ P(\xb) > 0 $ \\
$ \exists \ep_0 $ tel que $\ds \ab \in \N^N  \ \al_n \geq 1 \Rightarrow  \frac{{\p}^{\ab} P }{P}(\xb) \ll x_n^{-\ep_0} \ \ ( \xb \in J^N) $
\end{déf}

Notre premier résultat est de montrer que pour les polyn\^omes appartenant à cette classe, $ Z $ possède un prolongement {\bf holomorphe}:

\begin{theo}\label{existence du prolongement holomorphe de Z}
Soient $ Q,P_1,...,P_T \in \R[X_1,...,X_N] $. \\
On suppose que:

  $ \star \ \forall 1 \leq t \leq T \ P_t $ vérifie HDF 

 $ \ds \star \ \prod_{t=1}^T P_t(\xb) \xrightarrow[\substack{\xb \rightarrow + \infty \\ \xb \in J^N}]{} + \infty $. \\
Soit de plus $ \mub \in (\T \setminus \{ 1 \} )^N $. \\
Alors: \\
$ Z(Q,P_1,...,P_T,\mub,\cdot) $ possède un prolongement holomorphe à $ \C^T $.
\end{theo}

Le résultat suivant est crucial pour la suite: \\

\begin{theo}\label{lemme d'échange} { \bf  (lemme d'échange) } \\
Soient $ P_1,...,P_T,Q_1,...,Q_{T'} \in \R[X_1,...,X_N] $ et $ Q \in \R[X_1,...,X_N] $. \\
On suppose que:

   $ \star \ P_1,...,P_T,Q_1,...,Q_{T'} $ vérifient HDF

   $ \ds \star \ \prod_{t=1}^T P_t(\xb) \xrightarrow[ \substack{ |\xb| \rightarrow + \infty \\ \xb \in J^N}]{} + \infty $ \\

  $ \ds \star \ \prod_{t=1}^{T'} Q_t(\xb) \xrightarrow[ \substack{ |\xb| \rightarrow + \infty \\ \xb \in J^N}]{} + \infty $ \\
Soient de plus $ \mub \in (\T \setminus \{1\})^N $ et $ k_1,...,k_T,l_1,...,l_{T'} \in \N $. \\
Alors:
$$ Z \left(Q \prod_{t = 1}^{T'} {Q_t}^{l_t},P_1,...,P_T,\mub,-k_1,...,-k_T \right) = Z \left(Q \prod_{t = 1}^T {P_t}^{k_t},Q_1,...,Q_{T'},\mub,-l_1,...,-l_{T'} \right) $$
\end{theo}

Du lemme d'échange on déduit  des formules particulièrement simples pour les valeurs aux points de $ (- \N)^T $: \\

\begin{theo}\label{formule pour les valeurs de Z aux entiers négatifs}
Soient $ Q, P_1,...,P_T \in \R[X_1,...,X_N] $. \\
On suppose que: \\
a) $ P_1,...,P_T $ vérifient HDF \\
b) $ \ds \prod_{t=1}^T P_t(\xb) \xrightarrow[ \substack{ |\xb| \rightarrow + \infty \\ \xb \in J^N}]{} + \infty $. \\
Soient $ k_1,...,k_T \in \N $. On note $ \ds Q \prod_{t=1}^T P_t^{k_t}  = \sum_{\ab \in S} a_{\ab} \Xb^{\ab} $. \\
Soit de plus $ \mub \in (\T \setminus \{1\})^N $. \\
Alors:
$$ Z(Q,P_1,...,P_T,\mub,- k_1,...,- k_T) = \sum_{\ab \in S} a_{\ab} \prod_{n=1}^N \zeta_{\mu_n}(-\al_n) $$
\end{theo}

De ces formules on déduit les suivantes, qui permettent l'interpolation p-adique: \\

\begin{theo}\label{formule pour interpoler}
Soient $ Q,P_1,...,P_T \in \R[X_1,...,X_N] $. \\
On suppose que: \\
a) $ P_1,...,P_T $ vérifient $ HDF $, \\
b) $ \ds \prod_{t=1}^T P_t(\xb) \xrightarrow[\substack{ |\xb| \rightarrow +\infty \\ \xb \in J^N}]{} + \infty $. \\
Soit $ \mub \in ( \T \setminus \{ 1 \})^N $. \\
Alors pour tout $ \kb \in \N^N $ on a: 
$$ Z(Q,P_1,...,P_T,\mub,-\kb) = \frac{\mub^{\unb}}{(\unb - \mub)^{\unb}} \sum_{\lb \in \N^N}  \frac{1}{(\unb - \mub)^{\lb}}  \sum_{\jb \in \prod_{n=1}^N \{0,...,\ell_n \}} \left\{ (-1)^{|\jb|} \binom{\lb}{\jb} Q(-\jb) \prod_{t=1}^T P_t(-\jb)^{k_t} \right\} $$
formule dans laquelle la somme sur $ \lb $ est en fait une somme finie.
\end{theo}

On a en fin le:

\begin{theo}\label{interpolation}
Soit $ p $ un nombre premier.\\
On fixe un morphisme de corps de $ \C $ dans $ \C_p $, il sera sous entendu dans les écritures. \\
Soient $ Q,P_1,...,P_T \in \Z[X_1,...,X_N] $. \\
On suppose que: \\
a) $ P_1,...,P_T $ vérifient $ HDF $, \\
b) $ \ds \prod_{t=1}^T P_t(\xb) \xrightarrow[\substack{ |\xb| \rightarrow +\infty \\ \xb \in J^N}]{} + \infty $, \\
c) $ \forall t \in \{1,...,T \} \ \forall \jb \in \Z^N \ p \nmid P_t(\jb) $ . \\ 
Soit $ \mub \in ( \T \setminus \{ 1 \})^N $. \\
On suppose que $ \ds \forall n \in \{1,...,N \} \ |1-\mu_n|_p > p^{ - \frac{1}{p-1}} $. \\
Soit $ \rb \in \{0,..., p-1 \}^T $. \\
Alors il existe $ Z_p^{\rb}(Q,P_1,...,P_T, \mub, \cdot) \colon {\Z_p}^T \to \C_p $ continue telle que: \\
$ \forall \kb \in \N^T $ vérifiant $ \forall t \in \{1,...,T \} \ k_t = r_t \ [p-1] $, on ait: 
$$ Z_p^{\rb}(Q,P_1,...,P_T,\mub,-\kb) = Z(Q,P_1,...,P_T, \mub,-\kb). $$
\end{theo}

\part*{Quelques remarques}

Rappelons les définitions de deux classes usuelles de polyn\^omes:

\begin{notation}
Si $ P \in \R[X_1,...,X_N] $ s'écrit $ \ds P(\Xb) = \sum_{\ab \in \N^N} a_{\ab} \Xb^{\ab} $, alors on note $ \ds P^+(\Xb) = \sum_{\ab \in \N^N} |a_{\ab}| \Xb^{\ab} $
\end{notation}

\begin{déf}
$ P \in \R[X_1,...,X_N] \setminus \{0\} $ est dit non dégénéré si $ P(\xb) \asymp P^+(\xb) \ \ (\xb \in J^N) $ 
\end{déf}

Pour plus de détails sur la notion de polyn\^ome non dégénéré, on pourra consulter par exemple \cite{sargos 84}.

\begin{déf}
soit $ P \in \R[X_1,...,X_N] $. \\
P est dit hypoelliptique s'il vérifie les trois conditions suivantes: \\
$ P $  n'est pas constant \\
$ \forall \xb \in J^N \ P(\xb)>0 $ \\
$ \ds \forall \ab \in \N^N \setminus \{ \zerob \} \ \frac{{\p}^{\ab} P }{P}(\xb) \xrightarrow[\substack{|\xb| \rightarrow +\infty \\ \xb \in J^N}]{} 0  $
\end{déf}

\begin{rem}
Les polyn\^omes non dégénérés et les polyn\^omes hypoelliptiques vérifient HR. 
\end{rem}

  Dans \cite{essoua}, Essouabri a introduit une nouvelle classe de polyn\^omes qui contient les deux précédentes:

\begin{déf}
soit $ P \in \R[X_1,...,X_N] $ \\ 
On dit que $ P $ vérifie l'hypothèse $ H_0 S $ si: \\
$ \forall \xb \in J^N \ P(\xb) > 0 $ \\
$ \ds \forall \ab \in \N^N  \ \  \frac{{\p}^{\ab} P }{P}(\xb) \ll 1  \ \ ( \xb \in J^N) $
\end{déf}

  Clairement les méthodes \cite{essoua} de permettent de montrer le résultat suivant:

\begin{theo} 
Soient $ Q,P_1,...,P_T \in \R[X_1,...,X_N] $. \\
On suppose que:

  $ \star \ \forall 1 \leq t \leq T \ P_t $ vérifie $ H_0 S $. \\ 

 $ \ds \star \ \prod_{t=1}^T P_t(\xb) \xrightarrow[\substack{\xb \rightarrow + \infty \\ \xb \in J^N}]{} + \infty $. \\
Soit de plus $ \mub \in \T^N $. \\
Alors: \\
$ Z(Q,P_1,...,P_T,\mub,\cdot) $ possède un prolongement méromorphe à $ \C^T $.
\end{theo}

\begin{rem}
Il est clair que la classe des polyn\^omes vérifiant HDF  est incluse dans la classe des polyn\^omes vérifiant $ H_0 S $. L'exemple suivant montre que l'inclusion est stricte.
\end{rem}

\begin{ex}\label{H0S mais pas HR} 

Dans \cite{essoua}  Driss Essouabri remarque que $ P(X,Y) = (X-Y)^2 X + X \in \R[X,Y] $ vérifie $ H_0 S $ mais qu'il n'est pas hypoelliptique et qu'il est dégénéré. \\
En fait $ \ds \ \forall x \geq 1 \ \ \frac{\frac{\p P}{\p y}}{P}(x,x+1) = - 1 $ et donc $ P $ ne vérifie pas HDF.
\end{ex}

Ce polyn\^ome est intéressant car il montre que le théorème \ref{existence du prolongement holomorphe de Z} peut \^etre faux sous $ H_0 S $. C'est l'objet de l'exemple suivant:

\begin{ex}
Soit $ P(X,Y) = (X-Y)^2 X + X \in \R[X,Y] $. \\
Alors: \\
$ Z(1,P,-1,-1,\cdot) $ possède un prolongement méromorphe à $ \C $. \\
$ 1 $ est l'unique p\^ole du prolongement, il est simple de résidu $ \ds \frac{\pi}{\sinh(\pi)} $. \\
\end{ex}

\noindent { \bf Preuve}: \\
Durant cette preuve, on pose $ Z = Z(1,P,-1,-1,\cdot) $.
\begin{align*}
             Z(s) & = \sum_{m,n \geq 1} (-1)^m(-1)^n [(m-n)^2 m+m]^{-s} \\
                  & = \sum_{m,n \geq 1} (-1)^{m-n} m^{-s} [(m-n)^2 +1]^{-s} \\
                  & = \sum_{1 \leq m \leq n} (-1)^{m-n} m^{-s} [(m-n)^2 +1]^{-s} + \sum_{1 \leq n <m } (-1)^{m-n} m^{-s} [(m-n)^2 +1]^{-s}  
\end{align*}
En posant  $ n = m+u $ dans la première somme et $ m = n+u $ dans la deuxième, on obtient: 
\begin{align*}
             Z(s) & = \sum_{\substack{m \geq 1 \\ u \geq 0 }} (-1)^u m^{-s} (u^2 +1)^{-s} + \sum_{n,u \geq 1} (-1)^u (n+u)^{-s} (u^2+1)^{-s} \\
                  & = \zeta(s) \sum_{u \geq 0} (-1)^u (u^2+1)^{-s} + \sum_{u \geq 1} (-1)^u (u^2+1)^{-s} \sum_{n \geq 1} (n+u)^{-s} \\
                  & = \zeta(s) \sum_{u \geq 0} (-1)^u (u^2+1)^{-s} + \sum_{u \geq 1} (-1)^u (u^2+1)^{-s} \left[\zeta(s) - \sum_{1 \leq k \leq u} k^{-s} \right] \\ 
                  & = \zeta(s) \sum_{u \in \Z} (-1)^u (u^2+1)^{-s} - \sum_{1 \leq k \leq u} (-1)^u (u^2+1)^{-s} k^{-s} \\
                  & = \zeta(s) \sum_{u \in \Z} (-1)^u (u^2+1)^{-s} - \sum_{\substack{k \geq 1 \\ l \geq 0 }} (-1)^{k+l}[(k+l)^2+1]^{-s} k^{-s}
\end{align*}

\noindent Les deux observations suivantes permettent de conclure: \\
$ \ast $ c'est une application classique du théorème des résidus de montrer que $ \ds \sum_{u \in \Z} (-1)^u (u^2+1)^{-1} = \frac{\pi}{\sinh(\pi)} $, \\
$ \ast $ le théorème \ref{existence du prolongement holomorphe de Z} permet d'affirmer que $ \ds s \mapsto \sum_{\substack{k \geq 1 \\ l \geq 0 }} (-1)^{k+l}[(k+l)^2+1]^{-s} k^{-s} $ se prolonge holomorphiquement à $ \C $ . \\

\part*{Des lemmes sur les polyn\^omes à plusieurs variables}

\begin{lemme}\label{un polynome admet une dérivée partielle constante}
Soient $ P \in \R[X_1,...,X_N] $ et $ 1 \leq n \leq N $. \\
On suppose que $ P $ dépend effectivement de $ X_n $. \\
Alors:\\
 il existe $ \ab \in \N^N $ tel que $ \al_n \geq 1 $ et $ \p^{\ab}P $ soit constant et non nul. 
\end{lemme}

\noindent {\bf Preuve}: \\
On note S=supp(P) et $ \ds P(\Xb)= \sum_{\ab \in S } a_{\ab}{\Xb}^{\ab} $. \\
Notons $ d_n = deg_{X_n}P $, on a alors $ d_n \geq 1 $. \\
On prend $ \ab \in S $ tel que $ \al_n = d_n $ et $ |\ab| = \max \{|\bb| \ | \ \bb \in S, \ \be_n = d_n \} $ \\
Soit $ \bb \in S \setminus \{ \ab \} $,
\ \ si $ \be_n < d_n $ alors $ \p^{\ab} (\Xb^{\bb}) = 0 $,
\ \ si $ \be_n = d_n $ alors $ |\bb| \leq |\ab| $ et $ \bb \neq \ab $ donc $ \exists i \in [[1,N]] $ tel que $ \be_i < \al_i $, d'o\`u $ \p^\al(\Xb^{\bb}) = 0 $\\
On déduit de ceci que $ \p^{\ab} P = \p^{\ab} (a_{\ab} \Xb^{\ab}) = a_{\ab} \ab! $ donc $ \ab $ convient.

\begin{lemme}
 Si $ P \in \R[X_1,...,X_N] $ vérifie $ H_0 S $ alors $ P(\xb) \gg 1 \ \ ( \xb \in J^N) $.
\end{lemme}

\noindent {\bf Preuve}: \\
si $ P $ est constant, c'est clair. \\
On suppose $ P $ non constant; il existe alors $ n $ tel que $ P $ dépende effectivement de $ X_n $. \\
Le lemme \ref{un polynome admet une dérivée partielle constante} fournit $ \ab \in \N^N $ tel que $ \al_n \geq 1 $ et $ \p^{\ab} P $ soit constant et non nul. \\
$ P(\xb) \gg  \p^{\ab} P (\xb)  \ \ ( \xb \in J^N) $ donne alors le résultat. \\

\begin{lemme}\label{tarski saidenberg}

Soient $ 0 \leq N_1 \leq N $ et $ C $ un compact de $ \R^{N_1} $. \\
 Soit $ P \in \C(X_1,...,X_N) $. \\
On suppose que: \\
 $ R(\xb) \xrightarrow[\substack{ |\xb| \rightarrow + \infty \\ \xb \in C \times J^{N-N_1}}]{} 0 $. \\
Alors il existe $ \ep_0 > 0 $ tel que: \\
$ \ds  R(\xb) \ll \left( \prod_{n=N_1 + 1}^N x_n \right)^{-\ep_0} \ \ (\xb \in C \times J^{N-N_1}) $

\end{lemme}

\noindent {\bf Preuve}: \\
la preuve repose sur le principe de Tarski-Saidenberg qui est un outil classique de géométrie algébrique réelle. Pour plus de détails on pourra consulter par exemple \cite{essoua}.

\part*{Prolongement holomorphe des intégrales Y}

\begin{notation} 
pour $ 1 \leq t \leq T $ on note $ \eb_t = (0,...,0,1,0,...,0) \in \N^T $. \\
\end{notation}

\begin{déf}
Pour $ r \in \R $ on pose: \\
$ \br(r) = \{ f: [r,+\infty[ \to \C \ | \ \exists (f_n)_{n \in \N}  \ f_n: [r,+\infty[ \to \C \ C ^{\infty} $ bornée vérifiant $ f_0 = f \ f_{n+1}'=f_n \} $. \\
$ \br(r) $ est clairement un sous espace vectoriel de $ \C^{[r,+\infty[} $. \\
\end{déf}

\begin{lemme}\label{unicité de la suite associeé à une fonction de B(r)}
Soient $ r \in \R $ et $ f \in \br(r) $. \\
Alors: \\
1) il existe une unique suite $ (f_n)_{n \in \N} $ o\`u $ f_n \in \C^{[r,+\infty[} $ telle que: 

  $ \star \ \forall n \ f_n $ est $ C ^{\infty} $ bornée

  $ \star \ f_0 = f $

  $ \star \ \forall n \  f_{n+1}'=f_n $ \\
2) $ \forall n \in \N \ f_n \in \br(r) $.
\end{lemme}

\noindent {\bf Preuve}: \\ 
1) soient $  (f_n)_{n \in \N} $ et  $  (g_n)_{n \in \N} $ convenant. \\
Montrons par récurrence sur $ n \geq 0 $ que $ f_n = g_n $. \\
C'est clair pour $ n= 0 $. \\ 
Si l'on a  $ f_n = g_n $, alors: \\ 
$ f_{n+2}'' = f_{n+1}' = f_n = g_n =g_{n+1}' = g_{n+2}'' $ . \\
$ f_{n+2}'' = g_{n+2}'' $ donc $ f_{n+2} - g_{n+2} $ est une fonction affine sur $ [r,+\infty[ $, or elle est bornée, donc elle est constante, donc sa dérivée est nulle, c'est à dire $ f_{n+1} - g_{n+1} = 0 $, d'o\`u $ f_{n+1} = g_{n+1} $.  \\
2) $ f_n \in \br(r) $ est clair. \\ 

\noindent Le lemme suivant ne sera pas utilisé par la suite, mais répond à une question naturelle sur la classe $ \br(r) $. \\  

\begin{lemme} 
Soient $ r \in \R $ et $ f: [r,+\infty[ \to \C $. \\
Alors sont équivalents: \\
i) $ f \in \br(r) $ \\
ii) $ \forall n \ \exists g : [r,+\infty[ \to \C \ C ^{\infty} $ bornée telle que $ g^{(n)} = f. $ \\

\end{lemme}

\noindent {\bf Preuve}: \\
i) $ \Rightarrow $ ii) \\
Il suffit de remarquer que $ f_n^{(n)} = f $. \\
ii) $ \Rightarrow $ i) \\
pour tout $ n \in \N $ on choisit $ g_n: [r,+\infty[ \to \C  \ C ^{\infty} $ bornée telle que $ g_n^{(n)} =  f $. \\
$ (g_{n+1}')^{(n)} = f = g_n^{(n)} $ donc $ \exists P_n \in \C[X] $ de degré au plus $ n-1 $ tel que $ g_{n+1}' - g_n = P_n $. \\
 On note $ h_n = \Re (g_n) $ et $ R_n $ le polyn\^ome de $ \R [X] $ dont les coefficients sont les parties réelles de ceux de $ P_n $. Il vient $ h_{n+1}'- h_n = R_n $. \\
Supposons $ R_n $ non constant. \\
Si son coefficient dominant est strictement positif, alors $ R_n(x) \xrightarrow[ x \rightarrow + \infty]{} + \infty $ donc \\
 $ h_{n+1}'(x)  \xrightarrow[ x \rightarrow + \infty]{} + \infty $ puis $ h_{n+1}(x)  \xrightarrow[ x \rightarrow + \infty]{} + \infty $, ce qui est absurde puisque $ h_{n+1} $ est bornée. \\
On montre de m\^eme que le coefficient dominant de $ R_n $ ne peut \^etre strictement négatif. \\
On a une contradiction, donc $ R_n $ est constant. \\
En raisonnant de manière similaire sur les parties imaginaires on montre que $ \Im(g_{n+1}' - g_n) $ est constante. \\
On conclut de ce qui précède que $ g_{n+1}' - g_n $ est une fonction constante. \\
Pour tout $ n $ on pose $ f_n = g_{n+1}' $. \\
Alors: \\
$ \star \ f_n $ est  $ C^{\infty} $ bornée, \\
$ \star \ f_0 = g_1' = f $, \\
$ \star \ f_{n+1}' - f_n = g_{n+2}'' - g_{n+1}' = (g_{n+2}' - g_{n+1} )' = 0 $. \\
On en conclut que $ f \in \br(r) $. \\

\noindent Donnons deux exemples de familles de fonctions appartenant à $ \br(r) $, le premier est l'exemple "typique", le deuxième servira dans la preuve du théorème \ref{existence du prolongement holomorphe de Z}. \\

\begin{ex}\label{exemples de fonctions appartenant à B(r)}
Soit $ r \in \R $. \\
1) Soit $ f:[r,+\infty[ \to \C $ qui soit $ C ^{\infty} $ et périodique de valeur moyenne nulle. \\
Alors $ f \in \br(r) $. \\ 
2) Soient $ \al,\be \in \R $ et $ a \in \C $. \\
On suppose $ \ds \be \neq 0, \ \frac{\al}{\be} \notin \Z $ et $ |a| \neq 1 $.\\
Alors $ \ds f:[r,+\infty[ \to \C $ définie par $ \ds f(x) = \frac{\exp(i \al x)}{1-a\exp(i \be x)} $ est dans $ \br(r) $. \\
\end{ex}

\noindent {\bf Preuve}: \\
1) le  développement en série de Fourier de $ f $ donne le résultat. \\
2) $ \star $ cas $ |a| <1 $: \\
$ \ds f(x)=\exp( i \al x) \sum_{k=0}^{+ \infty} a^k \exp(i k \be x )=  \sum_{k=0}^{+ \infty} a^k \exp( i( \al + k \be) x) $ \\
On pose donc, pour $ \ds  n \in \N $, $ \ds f_n(x) =  \sum_{k=0}^{+ \infty} \frac{ a^k}{ (i( \al + k \be))^n }\exp( i( \al + k \be) x) $ \\
$ f_n $ est $ C^{\infty} $ bornée, $ f_{n+1}'=f_n \ f_0 = f $; donc $ f \in \br(r) $ . \\
$ \star $ cas $ |a| > 1 $: \\
$ \ds f(x) = \frac{a^{-1}  \exp(-i \be x ) \exp( i \al x) }{ a^{-1}  \exp(-i \be x ) - 1 } = -a^{-1} \frac{\exp(i( \al - \be) x )} {1- a^{-1}  \exp(i (- \be) x )} $ qui est dans $ \br(r) $ par le cas précédent. \\

\begin{theo}\label{existence du prolongement holomorphe de Y}
Soient $ Q,P_1,...,P_T \in \C[X_1,...,X_N] $ et $ 0 \leq N_1 \leq N $. \\
On suppose que: \\
a) $ \forall 1 \leq t \leq T $ on a: 

   $ \star \ \forall \xb \in [-1,1]^{N_1} \times J^{N-N_1} \ P_t(\xb) \notin \R_- $ 
 
  $ \star \ | P_t(\xb)| \gg 1 \ \ (\xb \in [-1,1]^{N_1} \times J^{N-N_1}) $ \\
b) $ \ds \prod_{t=1}^T |P_t(\xb)| \xrightarrow[\substack{ |\xb| \rightarrow +\infty \\ \xb \in [-1,1]^{N_1} \times J^{N-N_1}}]{} +\infty $ \\
c) $ \exists \ep > 0 $ tel que: $ \ds \ab \in \{ 0\} ^{N_1} \times \N^{N-N_1} \ \al_n \geq 1 \Rightarrow \frac{\p^{\ab} P_t}{P_t}(\xb) \ll x_n^{-\ep} \ \  ( \xb \in [-1,1]^{N_1} \times J^{N-N_1} ) $ \\
Soient de plus $ f: [-1,1]^{N_1} \rightarrow \C $ continue et $ f_{N_1+1},...,f_N \in \br(1) $. \\
On pose:
$$ \ds Y(Q,P_1,..,P_T, f_{N_1+1},..,f_N,f,\sb) = \int_{ [-1,1]^{N_1} \times J^{N-N_1} } Q(\xb)  \left( \prod_{t=1}^T P_t(\xb)^{-s_t} \right) f(x_1,..,x_{N_1}) \left( \prod_{n=N_1 + 1}^N f_n(x_n) \right) d\xb $$
Alors: \\
1) $ \exists \si_0 $ tel que: \\
$ \sb \mapsto  Y(Q,P_1,...,P_T, f_{N_1+1},...,f_N,f,\sb) $ existe et soit holomorphe sur $ \{ \sb \in \C^T \ | \ \forall 1 \leq t \leq T \ \si_t > \si_0 \} $ \\
2) $ Y(Q,P_1,...,P_T, f_{N_1+1},...,f_N,f,\cdot) $ possède un {\bf prolongement holomorphe à $ \C^T $}.
\end{theo}

\noindent {\bf Preuve} \\
Gr\^ace à \ref{tarski saidenberg} il existe $ \ep_0 > 0 $ tel que: \\
  $ \ds \prod_{t=1}^T |P_t(\xb)| \gg \left( \prod_{n=N_1 + 1}^N x_n \right)^{\ep_0} \ \ ( \xb \in [-1,1]^{N_1} \times J^{N-N_1} ) $ \\
Quitte à diminuer $ \ep_0 $ on peut bien sur imposer que l'on ait de plus: \\ 
 $ \ds \ab \in \{ 0\} ^{N_1} \times \N^{N-N_1} \ \al_n \geq 1 \Rightarrow \frac{\p^{\ab} P_t}{P_t}(\xb) \ll x_n^{-\ep_0} \ \  ( \xb \in [-1,1]^{N_1} \times J^{N-N_1} ) $ \\

\noindent 1) Preuve de l'existence de $ \si_0 $: \\ \\
Soit $ \si_0 \in \R $ que l'on va déterminer par la suite.  \\
Soit $ K $ compact de $ \C^T $ inclus dans $ \{ \sb \in \C^T \ | \ \forall 1 \leq t \leq T \ \si_t > \si_0 \} $. \\ 
$ \star $ Soit $ 1 \leq t \leq T $. \\
$ |P_t(\xb)| \gg 1 \ \ ( \xb \in [-1,1]^{N_1} \times J^{N-N_1}) $ donc $ \exists c>0 $ tel que $ \forall \xb \in [-1,1]^{N_1} \times J^{N-N_1} \ |P_t(\xb)| \geq c $. \\
$ \forall \xb \in [-1,1]^{N_1} \times J^{N-N_1} \ c^{-1}|P_t(\xb)| \geq 1 $ donc: $ \si_t > \si_0 \Rightarrow  (c^{-1}|P_t(\xb)|)^{\si_t} \geq  (c^{-1}|P_t(\xb)|)^{\si_0} $. \\
On en déduit $ |P_t(\xb)|^{\si_t} \gg |P_t(\xb)|^{\si_0}  \ \ ( \xb \in [-1,1]^{N_1} \times J^{N-N_1} \ \ \sb \in K) $ \\
$  |P_t(\xb)^{s_t}| =  |P_t(\xb)|^{\si_t} \exp[ - \tau_t \arg P_t(\xb)] $ donc $ |P_t(\xb)^{s_t}| \gg |P_t(\xb)|^{\si_t} \ \ ( \xb \in [-1,1]^{N_1} \times J^{N-N_1} \ \ \sb \in K) $. \\ 
On conclut de ce qui précède que: $ |P_t(\xb)^{-s_t}| \ll |P_t(\xb)|)^{- \si_0} $. \\
$ \star $ Il vient donc: $ \ds \prod_{t=1}^T P_t(\xb)^{-s_t} \ll \left[ \prod_{t=1}^T |P_t(\xb)| \right]^{-\si_0} \ \ ( \xb \in [-1,1]^{N_1} \times J^{N-N_1} \ \ \sb \in K) $ \\ 
On suppose désormais $ \si_0 > 0 $, alors $ \ds \prod_{t=1}^T P_t(\xb)^{-s_t} \ll \left( \prod_{n=N_1 + 1}^N x_n \right)^{- \si_0 \ep_0}  \ \ ( \xb \in [-1,1]^{N_1} \times J^{N-N_1} \ \ \sb \in K) $ \\ 
On note $ q = \max \{ \deg_{x_n}Q | \ N_1 + 1 \leq n \leq N \} $ (on peut évidemment supposer $ Q \neq 0 $). \\
On a alors:
$$  Q(\xb) \left( \prod_{t=1}^T P_t(\xb)^{-s_t} \right) f(x_1,..,x_{N_1}) \prod_{n=N_1 + 1}^N f_n(x_n)  \ll \left( \prod_{n=N_1 + 1}^N x_n \right)^{q - \si_0 \ep_0}  \ \ ( \xb \in [-1,1]^{N_1} \times J^{N-N_1} \ \ \sb \in K) $$ 
Ceci conduit à faire le choix suivant: $ \ds \si_0 = \frac{q+2}{\ep_0} > 0 $. \\
Le théorème garantissant l'holomorphie de fonctions définies à l'aide d'intégrales permet de conclure. \\

\noindent 2) Preuve de l'existence d'un prolongement holomorphe dans le cas $ N_1 = 0 $. \\

\noindent On adopte quelques conventions, valables durant la preuve de cette partie: \\
$ \star $ on dira qu'une fonction $ Y $ est combinaison entière des fonctions $ Y_1,...,Y_k $ s'il existe des fonctions $ \la,\la_1,...,\la_k: \C^T \rightarrow \C $ entières telles que $ Y = \la + \sum_{i=1}^k \la_i Y_i $. \\
$ \star $ les polyn\^omes $ P_1,...,P_T $ sont fixés pour toute la preuve, donc on abrège $ Y(Q,P_1,..,P_T, f_1,...,f_N,\cdot) $ en $ Y(Q, f_1,...,f_N,\cdot) $. \\
$ \star $ $ \br $ désigne $ \br(1) $ . \\

\noindent La preuve se fait par récurrence sur $ N $. \\

L'examen du passage du rang $ N-1 $ au rang $ N $ permet de montrer le résultat au rang $ N=1 $, le résultat au rang $ N=0 $ étant évident. Ceci dit, par commodité pour le lecteur, nous allons tout de m\^eme détailler la preuve au rang $ N=1 $. \\

\noindent Preuve du résultat au rang $ N=1 $. \\
Soient donc $ Q,P \in \C[X] $ o\`u $ P $ est non constant et vérifie $ \forall x \in [1, + \infty[ \ \ P(x) \notin \R_- $. \\
Soit de plus $ f \in \br $. \\
On veut montrer que $ Y(Q,f,\cdot) $ definie par $ \ds Y(Q,f,s) = \int_1^{+ \infty} Q(x)P(x)^{-s} f(x) dx $ se prolonge holomorphiquement à $ \C $. \\ 
On note $ p = \deg P $ et $ q = \deg Q $. \\
On constate que $ Y(Q,f,\cdot) $ est holomorphe sur $ \ds \left\{ s \in \C \ | \ \si > \frac{q+1}{p} \right\} $. \\
Montrons par récurrence sur $ m \geq 0 $ que $ \forall Q \in \C[X] \ \ \forall f \in \br \ \ Y(Q,f,\cdot) $ se prolonge holomorphiquement à $ \ds \left\{ s \in \C \ | \ \si > \frac{q+1-m}{p} \right\} $. \\
$ \star $ Au rang $ m = 0 $ le résultat est clair. \\
$ \star $ Supposons le résultat vrai au rang $ m $. \\
$ f \in \br $ donc le lemme \ref{unicité de la suite associeé à une fonction de B(r)} associe à $ f $ une suite de fonctions appartenant à $ \br $, on note $ f^1 $ le premier terme de cette suite. \\
Gr\^ace à une intégration par parties on a:
$$ Y(Q,f,s) = \left[ Q(x) P(x)^{-s} f^1(x) \right]_{x=1}^{x=+\infty} - \int_1^{+\infty} (Q'(x)P(x)^{-s} + Q(x) (-s) P'(x) P(x)^{-(s+1)}) f^1(x) dx $$
Ceci s'écrit: $ \ds Y(Q,f,s) = - Q(1)P(1)^{-s} f^1(1) - Y(Q',f^1,s) + s Y(QP',f^1,s+1) $. \\
Puisque $ \deg Q = q-1 $, par hypothèse de récurrence, $ Y(Q',f^1,\cdot) $ se prolonge holomorphiquement à $ \ds \left\{ s \in \C \ | \ \si > \frac{(q-1)+1-m}{p} \right\} $. \\
Puisque $ \deg(QP') = q+p-1 $, par hypothèse de récurrence, $ s \mapsto Y(QP',f^1,s+1) $ se prolonge holomorphiquement à  $ \ds \left\{ s \in \C \ | \ \si > \frac{(q+p-1)+1-m}{p} -1 \right\} $. \\ 
Or $ \ds \frac{(q-1)+1-m}{p} = \frac{(q+p-1)+1-m}{p} - 1 = \frac{q+1-(m+1)}{p} $, on a donc démontré le résultat au rang $ m+1 $. \\

\noindent Preuve du passage du rang $ N-1 $ au rang $ N $. \\ \\
Désormais on suppose le résultat vrai au rang $ N - 1 $ et l'on souhaite prouver le résultat au rang $ N $. \\ 
La preuve est découpée en 10 étapes. \\

\noindent Etape 1: \\
$ f_1 \in \br $ donc le lemme \ref{unicité de la suite associeé à une fonction de B(r)} associe à $ f_1 $ une suite de fonctions appartenant à $ \br $, on note $ f_1^1 $ le premier terme de cette suite. \\
On a alors que $ Y(Q,f_1,...,f_N,\sb) $ est  combinaison entière de $ \ds Y(\frac{\p Q}{\p x_1},f_1,...,f_N,\sb) $ et des \\
 $ \ds Y(Q \frac{\p P_t}{\p x_1},f_1^1,f_2,...,f_N,\sb + \eb_t) $ o\`u $ 1 \leq t \leq T $. \\

\noindent Preuve de l'étape 1:
\begin{align*}
Y(Q,f_1,...,f_N,\sb) & = \int_{J^N} Q(\xb)  \prod_{t=1}^T P_t(\xb)^{-s_t} \prod_{n=1}^N  f_n(x_n)  d\xb \\
                     & = \int_{J^{N-1}} \left\{ \int_1^{ + \infty} Q(\xb) \left( \prod_{t=1}^T P_t(\xb)^{-s_t} \right) f_1(x_1) dx_1 \right\}  \prod_{n=2}^N  f_n(x_n)  \prod_{n=2}^N dx_n 
\end{align*}
L'expression entre accolades est, gr\^ace à une intégration par parties par rapport à $ x_1 $, la différence de: $ \ds  \left[ Q(\xb) \left( \prod_{t=1}^T P_t(\xb)^{-s_t} \right) f_1^1(x_1) \right]_{x_1=1}^{x_1 = +\infty} $ \\
et de $ \ds \int_1^{ + \infty} \left( \frac{\p Q}{\p x_1}( \xb)  \prod_{t=1}^T P_t(\xb)^{-s_t} + Q(\xb) \sum_{t=1}^T (-s_t)\frac{\p P_t}{\p x_1}( \xb)  P_t ( \xb)^{-(s_t+1)} \prod_{r \neq t} P_r(\xb)^{-s_r} \right) f_1(x_1) dx_1 $ \\
On en déduit: 
\begin{align*}
Y(Q,f_1,...,f_N,\sb) = & -  \int_{J^{N-1}} Q(1,x_2,...,x_N) \left( \prod_{t=1}^T P_t(1,x_2,...,x_N)^{-s_t} \right) f_1^1(1) \prod_{n=2}^N f_n(x_n) \prod_{n=2}^N dx_n  \\
                      &  -  Y(\frac{\p Q}{\p x_1},f_1,...,f_N,\sb) + \sum_{t=1}^T s_t  Y(Q \frac{\p P_t}{\p x_1},f_1^1,f_2,...,f_N,\sb + \eb_t) 
\end{align*}
Les polyn\^omes de $ N-1 $ variables $ P_1(1,x_2,...,x_N),...,P_T(1,x_2,...,x_N) $ vérifient les hypothèses ad hoc, et donc, gr\^ace à l'hypothèse de récurrence, le terme défini par une intégrale sur $ J^{N-1} $ admet un prolongement holomorphe à $ \C^T $, ce qui permet de conclure. \\

\noindent Etape 2: \\
pour tout $ d \geq 1 \ \ Y(Q,f_1,...,f_N,\sb) $ est combinaison entière de $ \ds Y(\frac{\p^d Q}{\p x_1^d},f_1,...,f_N,\sb) $ et de fonctions du type: $ \ds Y(\frac{\p^i Q}{\p x_1^i} \frac{\p P_t}{\p x_1},g_1,...,g_N,\sb  + \eb_t) $ o\`u $  i \in \N, 1 \leq t \leq T $ et $ g_1,...,g_N \in \br $. \\

\noindent Preuve de l'étape 2: \\
La preuve se fait par récurrence sur $ d $. \\
Le rang $ d=1 $ résulte de l'étape 1. \\
Le passage de $ d $ à $ d+1 $ se fait en combinant le résultat au rang $ d $ et l'étape 1 appliquée au polyn\^ome $ \ds \frac{\p^d Q}{\p x_1^d} $. \\

\noindent Etape 3: \\

\noindent pour $ 1 \leq n \leq N \ \  Y(Q,f_1,...,f_N,\sb) $ est combinaison entière de fonctions du type: \\
$ \ds Y( \frac{ \p^i Q}{ \p x_n^i} \frac{ \p P_t}{\p x_n},g_1,...,g_N,\sb + \eb_t) $ o\`u $ i \in \N, 1 \leq t \leq T $ et $ g_1,...,g_N \in \br $. \\

\noindent Preuve de l'étape 3: \\
Il suffit bien  s\^ur de traiter le cas $ n = 1 $. \\
Pour obtenir le résultat pour $ n = 1 $ il suffit d'appliquer l'étape 2 avec $ d = \deg_{x_1} Q + 1 $. \\

\noindent Etape 4: \\
pour $ 1 \leq n \leq N, \ub \in \N^T $ et $ Q \in \C[X_1,...,X_N] $, on définit $ \er_{\ub}^n (Q) $ comme étant le sous espace vectoriel de $ \C[X_1,...,X_N] $ engendré par les polyn\^omes de la forme : $ \ds \p^{\bb} Q \prod_{k=1}^n \frac {\p^{|\ab_k|+1}P_{t_k}}{\p \xb^{\ab_k} \p x_k} $ o\`u: \\
 $ \bb \in \N^N, \ab_1,...,\ab_n \in \N^N $ et $ t_1,...,t_n \in [[1,T]] $ vérifient $ \forall 1 \leq t \leq T \ \ u_t = card \{1 \leq k \leq n | t_k=t \} $. \\
Il est clair que $ n \neq |\ub| \Rightarrow  \er_{\ub}^n (Q)= \{0 \} $. \\
On fait les deux observations suivantes:

  $ \star \ \er_{\ub}^n (Q) $ est stable par dérivation. 

  $ \star \ 1 \leq n \leq N-1, 1 \leq t \leq T $ et  $ \ds Q \in \C[\xb]  \Rightarrow \frac{\p P_t}{\p x_{n+1}}  \er_{\ub}^n (Q)  \subset  \er_{\ub+\eb_t}^{n+1} (Q) $. \\

\noindent Etape 5: \\
pour $ 1 \leq n \leq N $ et $ Q \in \C[X_1,...,X_N] $, $ Y(Q,f_1,...,f_N,\sb) $ est combinaison entière de fonctions du type: \\
$ \ds Y(R,g_1,..,g_N,\sb+\ub) $ o\`u $ \ub \in \N^T, R \in  \er_{\ub}^n (Q) $ et $ g_1,...,g_N \in \br $. \\

\noindent Preuve de l'étape 5: \\
la preuve se fait par récurrence sur $ n \in [[1,N]] $. \\
Pour $ n = 1 $ cela résulte de l'étape 3. \\
Supposons le résultat vrai au rang $ n $, o\`u $ n \in [[1,N-1]] $. \\
$ Y(Q,f_1,...,f_N,\sb) $ est donc combinaison entière de fonctions du type: \\
$ Y(R,g_1,...,g_N,\sb+\ub) $ o\`u $ \ub \in \N^T, R \in  \er_{\ub}^n (Q) $ et $ g_1,..,g_N \in \br $. \\
Par ailleurs, par l'étape 3, $ Y(R,g_1,...,g_N,\sb+\ub) $ est combinaison entière de fonctions du type: \\
 $ \ds Y( \frac{ \p^i R}{ \p x_{n+1}^i} \frac{ \p P_t}{\p x_{n+1}},h_1,...,h_N,\sb + \ub + \eb_t) $ o\`u $ i \in \N, 1 \leq t \leq T$ et $ h_1,...,h_N \in \br $. \\
Gr\^ace aux deux observations de l'étape 4 $ \ds \frac{ \p^i R}{ \p x_{n+1}^i} \frac{ \p P_t}{\p x_{n+1}} \in \er_{\ub+\eb_t}^{n+1} (Q) $, d'o\`u le résultat au rang $ n+1 $. \\

\noindent Etape 6: \\
pour $ \ub \in \N^T $ et $ Q \in \C[X_1,...,X_N] $, on note $ \er_{\ub} (Q) $ le sous espace vectoriel de $ \C[X_1,...,X_N] $ engendré par les polyn\^omes de la forme : $ \ds \p^{\bb} Q \prod_{t=1}^T \prod_{k \in F_t} \p^{f_t(k)} P_t $ o\`u: \\
$ \star \ \bb \in \N^N $ \\
$ \star $ les $ F_t $ sont des parties finies de $ \N $,  disjointes deux à deux, et vérifiant $ |F_t|=u_t $ \\
$ \star \ \forall 1\leq t \leq T \ \ f_t $ est une fonction de $ F_t $ dans $ \N^N $ \\
$ \star $ on peut associer aux $ f_t $ des parties finies de $ \N $, $ D_1,...,D_N $ disjointes deux à deux et telles que:

  $ \ast \ |D_1|=...=|D_N| $, 

  $ \ds \ast \ \bigsqcup_{n=1}^N D_n = \bigsqcup_{t=1}^T F_t $ 

  $ \ds \ast \ 1 \leq t \leq T, \ 1 \leq n \leq N $ et $  k \in D_n \cap F_t \Rightarrow f_t(k) \in \N^{n-1} \times \N^* \times \N^{N-n} $. \\

\noindent Remarquons que $ \er_{\ub} (Q) $ est stable par dérivation. \\

\noindent Etape 7: \\
$ R \in \er_{\ub} (Q) $ et $ S \in \er_{\vb} (R) \Rightarrow S \in \er_{\ub+ \vb} (Q) $. \\ \\
Preuve de l'étape 7: \\
$ S $ est une combinaison linéaire de termes de la forme: $ \ds  \p^{\bb} R \prod_{t=1}^T \prod_{k \in F'_t} \p^{f'_t(k)} P_t $ o\`u: \\
$ \bb \in \N^N, |F'_t|=v_t, f'_t: F'_t \rightarrow \N^N $ et $ D'_1,...,D'_N $ sont comme à l'étape 6. \\
$ R \in  \er_{\ub} (Q) $ donc $ \p^{\bb} R \in \er_{\ub} (Q) $ donc $ \p^{\bb} R $ est combinaison linéaire de termes de la forme: $ \ds  \p^{\gb} Q \prod_{t=1}^T \prod_{k \in F_t} \p^{f_t(k)} P_t $ \\
o\`u $ \gb \in \N^N, |F_t|=u_t, f_t: F_t \rightarrow \N^N $ et $ D_1,...,D_N $ sont comme à l'étape 6. \\
On peut imposer que $ \forall t,t' \ F_t \cap F'_{t'} = \emptyset $. Ceci entraine: \\
$ \forall n,t \ D_n \cap F'_t = F_t \cap D'_n = \emptyset $ et $ \forall n,n' \ D_n \cap D'_{n'} = \emptyset $. \\
Pour conclure il nous suffit de voir que: \\
 $ \ds U \overset{\text{déf}}{=} \p^{\gb} Q \left( \prod_{t=1}^T \prod_{k \in F_t} \p^{f_t(k)} P_t \right) \left( \prod_{t=1}^T \prod_{k \in F'_t} \p^{f'_t(k)} P_t \right) $ est dans $ \er_{\ub+ \vb} (Q) $. \\
Pour $ 1 \leq t \leq T $ on définit $ g_t: F_t \sqcup F'_t \rightarrow \N^N $ par: $ g_t(k)=f_t(k) $ si $ k \in F_t $ et $ g_t(k)=f'_t(k) $ si $ k \in F'_t $. \\
Il vient alors que: $ \ds U = \p^{\gb} Q \prod_{t=1}^T \prod_{k \in F_t\sqcup F'_t} \p^{g_t(k)} P_t $. \\
Sous cette forme on va voir que $ U \in \er_{\ub+ \vb} (Q) $. \\ 
$ \star $ Les $ F_t \sqcup F'_t $ sont disjointes deux à deux et $ \forall \ 1 \leq t \leq T \ \ |F_t \sqcup F'_t|=u_t+v_t $. \\
$ \star $ Les $ D_n \sqcup D'_n $ sont disjointes deux à deux, $ \ds \bigsqcup_{t=1}^T ( F_t\sqcup F'_t) = \bigsqcup_{n=1}^N ( D_n \sqcup D'_n) $ et $ |D_1 \sqcup D'_1| =...=|D_N \sqcup D'_N| $ \\
$ \star $ Si $ k \in ( D_n \sqcup D'_n) \cap ( F_t \sqcup F'_t) = (D_n \cap F_t) \sqcup  (D'_n \cap F'_t) $, alors: 

  $ \ast $ soit $ k \in D_n \cap F_t $ et alors $ g_t(k)=f_t(k) \in \N^{n-1} \times \N^* \times \N^{N-n} $, 

  $ \ast $ soit $ k \in D'_n \cap F'_t $ et alors $ g_t(k)=f'_t(k) \in \N^{n-1} \times \N^* \times \N^{N-n} $. \\
On en conclut que l'on a bien $ U \in \er_{\ub+ \vb} (Q) $. \\

\noindent Etape 8: \\
$ Q \in \C[X_1,...,X_N] $ et $ \ub \in \N^T \Rightarrow \er_{\ub}^N(Q) \subset \er_{\ub}(Q) $. \\

\noindent Preuve de l'étape 8: \\
On pose $ \ds S = \p^{\bb} Q \prod_{k=1}^N \frac {\p^{|\ab_k|+1} P_{t_k}}{\p \xb^{\ab_k} \p x_k} $ o\`u: \\
$ \bb \in \N^N, \ab_1,...,\ab_N \in \N^N $ et $ t_1,...,t_N \in [[1,T]] $ vérifient $ \forall 1 \leq t \leq T \ \ u_t = card \{1 \leq k \leq N | t_k = t \} $. \\
Pour conclure il suffit de montrer que $ S \in \er_{\ub}(Q) $. \\
Pour $ 1 \leq t \leq T $ on pose $ F_t = \{ 1 \leq k \leq N \ | \ t_k = t \} $; on a alors que $ |F_t| = u_t $. \\
On constate que les $ F_t $ sont deux à deux disjoints et que $ \ds \bigsqcup_{t=1}^T F_t = [[1,N]] $. \\
On définit $ f_t : F_t \to \N^N $ par $ f_t(k) = \ab_k+ \eb_k $. \\
On pose $ D_n = \{ n\} $. \\
On constate alors :

  $ \ast \ |D_1| = ... = |D_N| $ \\

  $ \ds \ast \ \bigsqcup_{n=1}^N D_n = [[1,N]] $ \\

  $ \ds \ast $ si $ k \in D_n \cap F_t $ alors $ k = n $ et donc $ f_t(k) = \ab_n + \eb_n \in \N^{n-1} \times \N^* \times \N^{N-n} $. \\

\noindent $ \ds S = \p^{\bb} Q \prod_{t=1}^T \prod_{ k \in F_t} \frac {\p^{|\ab_k|+1} P_{t_k}}{\p \xb^{\ab_k} \p x_k} = \p^{\bb} Q \prod_{t=1}^T \prod_{ k \in F_t} \p^{f_t(k)} P_t, $ donc $ S \in \er_{\ub}(Q) $. \\

\noindent Etape 9: \\
\\
Soient $ m \geq 1 $ et $ Q \in \C[X_1,...,X_N] $. \\
Alors $ Y(Q,f_1,...,f_N,\sb) $ est combinaison entière de fonctions du type:\\
 $ Y(R,g_1,...,g_N, \sb + \ub) $ o\`u $ \ub \in \N^T,|\ub| = mN, R \in \er_{\ub}(Q) $ et $ g_1,...,g_N \in \br $. \\

\noindent Preuve de l'étape 9: \\
\\
on raisonne par récurrence sur $ m \geq 1 $. \\
$ \ast $ pour $ m=1 $: \\
par l'étape 5 \ $ Y(Q,f_1,...,f_N,\sb) $ est combinaison entière de fonctions du type $ Y(R,g_1,...,g_N, \sb + \ub) $ o\`u $ \ub \in \N^T, R \in \er_{\ub}^N(Q) $ et $ g_1,...,g_N \in \br $. \\
On peut supposer $ |\ub| = N $ (car si $ |\ub| \neq N $ alors $ \er_{\ub}^N(Q) = \{0\} $). L'étape 8 donne donc le résultat. \\
$ \ast $ supposons le résultat vrai au rang $ m \geq 1 $. \\
$ Y(Q,f_1,...,f_N,\sb) $ est combinaison entière de fonctions du type $ Y(R,g_1,...,g_N, \sb + \ub) $ o\`u: \\
$ \ub \in \N^T, |\ub| = mN, R \in \er_{\ub}(Q) $ et $ g_1,...,g_N \in \br $. \\
Par ailleurs, le résultat pour $ m=1 $ donne que $ Y(R,g_1,...,g_N, \sb + \ub) $ est combinaison entière de fonctions du type $  Y(S,h_1,...,h_N, \sb + \ub + \vb) $ o\`u $ \vb \in \N^T, |\vb| = N, S \in \er_{\vb}(R) $ et $ h_1,...,h_N \in \br $. \\
L'étape 7 donne alors $ S \in \er_{\ub + \vb}(Q) $, mais $ |\ub+\vb| = (m+1)N $, d'o\`u le résultat au rang $ m+1 $. \\

\noindent Etape 10: conclusion \\
\\
On fixe $ Q \in \C[X_1,...,X_N] $ et $ f_1,...,f_N \in \br $ jusqu'à la fin. \\
Soit $ m \geq 1 $. \\
Par l'étape 9 $ Y(Q,f_1,...,f_N,\sb) $ est combinaison entière de fonctions du type:\\
 $ Y(R,g_1,...,g_N, \sb + \ub) $ o\`u $ \ub \in \N^T,|\ub| = mN, R \in \er_{\ub}(Q) $ et $ g_1,...,g_N \in \br $. \\
$ R \in \er_{\ub}(Q) $ donc $ R $ s'écrit: comme une combinaison linéaire de polyn\^omes de la forme: \\
$ \ds \p^{\bb} Q \prod_{t=1}^T \prod_{k \in F_t} \p^{f_t(k)} P_t $ avec $ \bb \in \N^N, |F_t| = u_t, f_t :F_t \to \N^N $ et $ D_1,...,D_N $ comme à l'étape 6. \\
On a alors $ \forall n \ |D_n| = m $. \\ 
Il vient:
\begin{align*}
\prod_{t=1}^T \prod_{k \in F_t} \p^{f_t(k)} P_t(\xb) & = \prod_{t=1}^T \prod_{n=1}^N \prod_{k \in F_t \cap D_n} \p^{f_t(k)} P_t(\xb) \\
                                                     & \ll \prod_{t=1}^T \prod_{n=1}^N \prod_{k \in F_t \cap D_n} x_n^{-\ep_0} P_t(\xb) \ \ (\xb \in J^N) \\
                                                     & \ll \prod_{t=1}^T \prod_{n=1}^N  (x_n^{-\ep_0} P_t(\xb))^{|F_t \cap D_n|} \ \ (\xb \in J^N) \\
                                                     & \ll \prod_{n=1}^N \prod_{t=1}^T x_n^{-\ep_0 |F_t \cap D_n|} \prod_{t=1}^T \prod_{n=1}^N P_t(\xb)^{|F_t \cap D_n|}  \ \ (\xb \in J^N) \\ 
                                                     & \ll \prod_{n=1}^N  x_n^{-\ep_0 |D_n|} \prod_{t=1}^T P_t(\xb)^{|F_t|}  \ \ (\xb \in J^N) \\
                                                     & \ll \prod_{n=1}^N  x_n^{-\ep_0 m} \prod_{t=1}^T P_t(\xb)^{u_t} \ \ (\xb \in J^N) 
\end{align*}

\noindent On pose $ q = \max \{ \deg_{X_n} Q | \ 1 \leq n \leq N \} $ (on peut évidemment supposer $ Q \neq 0 $). \\
On pose $ p = \max \{ \deg_{X_n} P_t | \ 1 \leq n \leq N \ 1 \leq t \leq T\} $. \\
Soit $ a \in \R $ que l'on va déterminer par la suite. \\
Soit $ K $ compact de $ \C^T $ inclus dans $ \{ \sb \in \C^T \ | \ \forall 1 \leq t \leq T \ \si_t > -a \} $. \\ 
$ \star $ Soit $ 1 \leq t \leq T $. \\
Comme lors de la preuve de l'existence de $ \si_0 $ on montre $ |P_t(\xb)^{-\si_t}| \ll | P_t(\xb)|^a  \ \ ( \xb \in J^N \ \ \sb \in K) $. \\ 
On suppose désormais $ a > 0 $, alors $ \ds P_t(\xb)^a \ll \left( \prod_{n=1}^N x_n \right)^{pa} \ \ ( \xb \in J^N ) $. \\ 
Des inégalités précédentes on déduit: $ \ds P_t(\xb)^{-s_t} \ll \left( \prod_{n=1}^N x_n \right)^{pa} \ \ ( \xb \in J^N \ \ \sb \in K) $. \\
Posons $ \ds S = \p^{\bb} Q \prod_{t=1}^T \prod_{k \in F_t} \p^{f_t(k)} P_t $; en combinant ce qui précède, il vient alors:
\begin{align*}
  S(\xb) \prod_{t=1}^T P_t(\xb)^{-(s_t + u_t)} & \ll \p^{\bb} Q(\xb) \prod_{n=1}^N  x_n^{-\ep_0 m} \prod_{t=1}^T P_t(\xb)^{u_t} \prod_{t=1}^T P_t(\xb)^{-(s_t + u_t)} \ \ (\xb \in J^N \ \ \sb \in K)  \\
                                               & \ll \left( \prod_{n=1}^N x_n \right)^q \left( \prod_{n=1}^N  x_n \right)^{-\ep_0 m} \left( \prod_{n=1}^N x_n \right)^{Tpa} \ \ ( \xb \in J^N \ \ \sb \in K) \\
                                               & \ll \left( \prod_{n=1}^N x_n \right)^{q + Tpa - \ep_0 m} \ \ ( \xb \in J^N \ \ \sb \in K) 
\end{align*}

\noindent On suppose désormais que $ \ds m > \frac{q+2}{\ep_0} $. \\
On choisit $ \ds a = \frac{ \ep_0 m - (q +2)}{Tp} $; ceci est bien strictement positif. \\
Ce qui précède montre que $ Y(S,g_1,...,g_N, \sb + \ub) $ est holomorphe sur $ \ds \{ \sb \in \C^T \ | \ \forall 1 \leq t \leq T \ \si_t > - a \} $. \\
On en déduit que: \\
$ Y(Q,f_1,...,f_N,\cdot) $ possède un prolongement holomorphe à $ \ds \left\{ \sb \in \C^T \ | \ \forall 1 \leq t \leq T \ \si_t > \frac{ q+2 - \ep_0 m }{Tp} \right\} $.
Ceci étant vrai pour tout $ \ds m > \frac{2}{\ep_0} $, $ Y(Q,f_1,...,f_N,\cdot) $ possède un prolongement holomorphe à $ \C^T $. \\

\noindent Comme le montre l'exemple suivant, le théorème \ref{existence du prolongement holomorphe de Y} peut ne plus \^etre vrai si l'on supprime l'hypothèse c).

\begin{ex}
On reprend le polyn\^ome de l'exemple \ref{H0S mais pas HR}: $ P(X,Y) = (X-Y)^2 X + X $. \\
On définit $ f_1: J \to \C $ par $ f_1 (x) = e^{i x} $ et $ f_2: J \to \C $ par $ f_1(y) = e^{-i y} $. $ f_1 $ et $ f_2 $ appartiennent à $ \br(1) $. \\
Alors: \\
$ Y(1,P,f_1,f_2,\cdot) $ possède un prolongement méromorphe à $ \C $. \\
$ 1 $ est l'unique p\^ole du prolongement, il est simple de résidu égal à $ \ds \frac{\pi}{e} $. \\
\end{ex}

\noindent {\bf Preuve}: \\
Par définition $ \ds Y(1,P,f_1,f_2,s) = \int_{J^2} P(x,y)^{-s} e^{i(x-y)} \ dxdy $ . \\
On pose $ \ds Y_1(s) =  \int_{ \{(x,y)|1<x<y \} } P(x,y)^{-s} e^{i(x-y)}  dxdy $. \\
Soit $ f: ]1, + \infty[ \times \R_+^* \to  \{ (x,y)|1<x<y \} $ définie par $ f(u,v) = (u,u+v) $. \\
$ f $ est un $ C^1 $ difféomorphisme dont le jacobien vaut partout $ 1 $. \\
En utilisant $ f $, on voit que: $ \ds Y_1(s) = \int_{ ]1, + \infty[ \times \R_+^* } [(u-(u+v))^2 u+ u]^{-s} e^{i(u-(u+v))} dudv $ \\
Donc $ \ds Y_1(s) = \int_1^{+ \infty} u^{-s} du  \int_0^{+ \infty}  (v^2+1)^{-s} e^{-iv} dv = \frac{1}{s-1} \int_0^{+ \infty} (v^2+1)^{-s} e^{-iv} dv $ \\
On pose $ \ds Y_2(s) =  \int_{ \{(x,y) | 1<y<x \} } P(x,y)^{-s} e^{i(x-y)}  dxdy $. \\  
Soit $ g:]1, +\infty[ \times \R_+^* \to  \{(x,y)|1<y<x \} $ définie par $ g(u,v) = (u+v,u) $. \\
$ g $ est un $ C^1 $ difféomorphisme dont le jacobien vaut partout $ -1 $. \\ 
En utilisant $ g $, on voit que: $ \ds Y_2(s) = \int_{ ]1, + \infty[ \times \R_+^* } [(u+v-u)^2 (u+v)+(u+v)]^{-s}e^{i(u+v-u)} \ dudv $ 
\begin{align*}
 \text{d'o\`u:} \ Y_2(s) & = \int_{ ]1, + \infty[ \times \R_+^* } (v^2+1)^{-s} (u+v)^{-s} e^{iv} dudv  \\
                         & = \int_0^{+ \infty } (v^2+1)^{-s} e^{iv} \left\{ \int_1^{+ \infty} (u+v)^{-s} du \right\} dv  \\
                         & = \int_0^{+ \infty }  (v^2+1)^{-s} e^{iv}   \frac{(1+v)^{-s+1}}{s-1}  dv \\
                         & = \frac{1}{s-1} \int_0^{+ \infty }  (v^2+1)^{-s} (v+1)^{-s+1} e^{iv} dv 
\end{align*}
Posons $ \ds Y(s) = \int_0^{+ \infty} (v^2+1)^{-s} e^{-iv} dv + \int_0^{+ \infty }  (v^2+1)^{-s} (v+1)^{-s+1} e^{iv} dv $. \\
Le théorème \ref{existence du prolongement holomorphe de Y} permet d'affirmer que $ Y $ admet un prolongement holomorphe à $ \C $. \\
On a $ \ds Y(1,P,f_1,f_2,s) = \frac{1}{s-1} Y(s) $, on va donc chercher à évaluer $ Y(1) $. \\
$ \ds Y(1) = \int_0^{+ \infty} (v^2+1)^{-1} e^{-iv} dv + \int_0^{+ \infty }  (v^2+1)^{-1} e^{iv} dv = \int_{- \infty}^{+ \infty }  (v^2+1)^{-1} e^{iv} dv $ \\
C'est une application classique du théorème des résidus de montrer que cette dernière intégrale vaut $ \ds \frac{\pi}{e} $, d'o\`u le résultat. \\

\part*{Domaine de convergence de Z}

\begin{notation}
Pour $ S \subset \R^T, int(S) $ désigne l'intérieur (dans $ \R^T $) de $ S $. 
\end{notation}

\begin{déf}\label{domaine de convergence de Z}
Soient $ Q,P_1,...,P_T \in \R[X_1,...,X_N] $ tels que $ \forall 1 \leq t \leq T \ \ \forall \xb \in J^N \ P_t(\xb) > 0 $. \\
On pose: \\
$ \Cr(Q,P_1,...,P_T) = \{ (\si_1,...,\si_T) \in \R^T \ | \ Z(Q,P_1,...,P_T,\unb,\si_1,...,\si_T) \ \mbox{converge} \} $ . \\

\end{déf}

\begin{rem}
Si, de plus, $ \mub $ appartient à $ {\T}^N $, alors on a:\\
$ Z(Q,P_1,...,P_T,\mub,s_1,...,s_T) $ converge $ \Leftrightarrow (\si_1,...,\si_T) \in \Cr(Q,P_1,...,P_T) $. 
\end{rem}

\begin{prop}\label{convexité du domaine de convergence de Z} 
Soient $ Q,P_1,...,P_T \in \R[X_1,...,X_N] $ tels que $ \forall 1 \leq t \leq T \ \ \forall \xb \in J^N \ P_t(\xb) > 0 $. \\
Alors $ \Cr(Q,P_1,...,P_T) $ est convexe (et donc connexe). \\

\end{prop}

\noindent {\bf Preuve}: \\
Soient $ \sib, \sib' \in \Cr(Q,P_1,...,P_T) $ et $ \la \in [0,1] $. \\
Fixons $ \mb \in \N^{*N} $ et posons $ a_t = P_t(\mb)^{-1} $, il vient alors:
$$ \prod_{t=1}^T P_t(\mb)^{- (\la \si_t + (1 - \la) \si_t')} = \prod_{t=1}^T a_t^{ \la \si_t + (1 - \la) \si_t'}  = \left( \prod_{t=1}^T a_t^{\si_t} \right)^{\la} \left( \prod_{t=1}^T a_t^{\si_t'} \right)^{1 - \la} $$
En utilisant l'inégalité $ a^{\la} b^{1-\la} \leq \la a + (1-\la)b $, valable pour $ a,b >0 $ et $ \la \in [0,1] $, on voit que:
$$ \left( \prod_{t=1}^T a_t^{\si_t} \right)^{\la} \left( \prod_{t=1}^T a_t^{\si_t'} \right)^{1 - \la}\leq \la \prod_{t=1}^T a_t^{\si_t} + (1 - \la) \prod_{t=1}^T a_t^{\si_t'} $$
On conclut de ce qui précède que $ \la \sib + (1 - \la) \sib' \in \Cr(Q,P_1,...,P_T) $.

\begin{lemme}\label{une propriété du domaine de convergence de Z}
Soient $ Q,P_1,...,P_T \in \R[X_1,...,X_N] $ tels que $ \forall 1 \leq t \leq T \ P_t(\xb) \gg 1 \ \ (\xb \in J^N) $. \\
Soit $ \ub \in \R_+^T $. \\
Alors: $ \Cr(Q,P_1,...,P_T) + \ub \subset \Cr(Q,P_1,...,P_T) $. 
\end{lemme}

\noindent {\bf Preuve}: \\
c'est clair. \\

\begin{cor}\label{une propriété de l'intérieur du domaine de convergence de Z}
Soient $ Q,P_1,...,P_T \in \R[X_1,...,X_N] $ tels que $ \forall 1 \leq t \leq T \ P_t(\xb) \gg 1 \ \ ( \xb \in J^N) $. \\
Soit $ \ub \in \R_+^{*T} $. \\
Alors: $ \Cr(Q,P_1,...,P_T) + \ub \subset int(\Cr(Q,P_1,...,P_T)) $. 
\end{cor}

\noindent {\bf Preuve}:
cela découle du lemme \ref{une propriété du domaine de convergence de Z}.

\begin{prop}\label{pour sigma grand on est dans le domaine de convergence de Z}
Soient $ Q,P_1,...,P_T \in \R[X_1,...,X_N] $ tels que $ \forall 1 \leq t \leq T \ P_t(\xb) \gg 1 \ \ ( \xb \in J^N) $. \\
Soit $ 1 \leq T_0 \leq T $. \\
On suppose que  $ \ds \prod_{t=1}^{T_0} P_t(\xb) \xrightarrow[\substack{\xb \rightarrow + \infty \\ \xb \in J^N}]{} + \infty $. \\
Soit $ \si_{T_0 + 1},...,\si_T \in \R $. \\
Alors: \\
il existe $ \si_0 \in \R $ tel que: $ \si_1,...,\si_{T_0} \geq \si_0 \Rightarrow (\si_0,...,\si_{T_0},\si_{T_0 + 1},...,\si_T) \in int(\Cr(Q,P_1,...,P_T)) $. \\
\end{prop}

\noindent {\bf Preuve}: \\
Soit $ \si_0 \in \R $ tel que $ \ds \prod_{t=1}^{T_0} P_t(\xb)^{\si_0 - 1} \prod_{t = T_0 + 1 }^T P_t(\xb)^{\si_t - 1} \gg  \left( \prod_{n=1}^N x_n \right)^2 \ \ (\xb \in J^N) $. \\
Alors $ (\si_0,...,\si_0,\si_{T_0 + 1},...,\si_T) - \unb \in \Cr(Q,P_1,...,P_T) $. \\
Le corollaire \ref{une propriété de l'intérieur du domaine de convergence de Z} permet alors de conclure. \\

\begin{prop}\label{holomorphie de Z}
Soient $ Q,P_1,...,P_T \in \R[X_1,...,X_N] $ tels que $ \forall 1 \leq t \leq T \ P_t(\xb) \gg 1 \ \ ( \xb \in J^N) $. \\
Soit $ \mub \in {\T}^N $. \\
Alors $ Z(Q,P_1,...,P_T,\mub,\cdot) $ est holomorphe sur $ int(\Cr(Q,P_1,...,P_T)) + i \R^T $. 
\end{prop}

\noindent Preuve: \\
a) Montrons que pour $ C $ compact inclus dans $ int(\Cr(Q,P_1,...,P_T)) $, $ Z(Q,P_1,...,P_T,\mub,\cdot) $ converge normalement sur $ C $. \\
Soit $ \sib \in C $; alors il existe $ \ep >0 $ tel que $ \sib - \ep \unb \in \Cr(Q,P_1,...,P_T) $. \\
Par hypothèse $ \exists c \in ]0,1[ $ tel que $ \forall 1 \leq t \leq T \ \ \forall \xb \in J^N\ P_t(\xb) \geq c $. \\
Soit $ 1 \leq t \leq T $. \\
Pour $ \ub \in ]0,2\ep[^N $ et $ \mb \in \N^{*N} $ on a: \\ 
$ P_t(\mb)^{-u_t} \leq c^{-u_t} $ et $ c^{-u_t} \leq c ^{-2 \ep} $; et donc $ P_t(\mb)^{-u_t} \leq c ^{-2 \ep} $. \\
On déduit de ce qui précède que: 
$$ |Q(\mb)| \prod_{t=1}^T P_t(\mb)^{-(\si_t - \ep + u_t)} \ll |Q(\mb)| \prod_{t=1}^T P_t(\mb)^{-(\si_t - \ep)}  \ \ (\mb \in \N^{*N}, \ub \in ]0,2 \ep[^N) $$
Ceci implique que $ Z(Q,P_1,...,P_T,\mub,\cdot) $ converge normalement sur $ \sib - \ep \unb + ]0,2\ep[^N $. \\
Pour tout $ \sib \in C $ on a trouvé un ouvert de $ \R^N $ contenant $ \sib $ et sur lequel il y a convergence normale; $ C $ étant compact, il y a convergence normale sur $ C $. \\
b) Conclusion. \\
On définit $ p: \C^N \to \R^N $ par $ p(\sib + i \taub) = \sib $. \\
Soit $ K $ un compact de $ int(\Cr(Q,P_1,...,P_T)) + i \R^T $. \\
$ p(K) $ est alors un compact de $ int(\Cr(Q,P_1,...,P_T)) $; par a) il y a convergence normale sur $ p(K) $, on en déduit la convergence normale sur $ K $. \\
De la convergence sur tout compact inclus dans $ int(\Cr(Q,P_1,...,P_T)) + i \R^T $ on déduit l'holomorphie sur $ int(\Cr(Q,P_1,...,P_T)) + i \R^T $.

\part*{Représentation intégrale}

\begin{notation} 
On pose $ \zerob = (0,...,0) \in \N^T $. \\
Si $ \ga: [a,b] \to \C $, on définit $ \ga^{-} : [a,b] \to \C $ par $ \ga^{-}(x) = \ga(a+b-x) $. \\
\end{notation}

\begin{lemme}\label{formule de représentation intégrale en dimension 1}

Pour $ \ep >0 $ on définit: \\
$ \ds \la_\ep: \left[\frac{3}{2}, + \infty \right[  \to \C $ par $ \la_\ep(x) = x + i \ep $ et $ \ds \om_\ep: \left[\frac{3}{2}, + \infty \right[  \to  \C $ par $ \om_\ep(x) = x - i \ep $. \\
Soit $ k \in \N^* $. \\
On note: $ \ds \la_{\ep,k}= \la_{\ep | [ \frac{3}{2}, k + \frac{1}{2}]} $ et $ \ds \om_{ \ep,k} = \om_{\ep | [ \frac{3}{2}, k + \frac{1}{2}]} $. \\ 
On définit $ \ds \ga_{ \ep,k}: [-1,1] \to \C $ par $ \ds \ga_{ \ep,k}(x) = k + \frac{1}{2} + i \ep x $. \\
On note $ \ds R_{\ep , k} = \left[ \frac{3}{2}, k + \frac{1}{2} \right] + i[ - \ep, \ep ] $. \\
On définit $ e : \C \to \C $ par $ e(z) = \exp(2 i \pi z) $ . \\
Soit $ f: U \rightarrow \C $ holomorphe o\`u $ U $ est un ouvert simplement connexe de $ \C $ contenant  $ R_{\ep , k} $ . \\
Alors: $ \ds  \sum_{m=2}^k f(m) = \int_{ \ga_{ \ep,1}^-} \frac{ f(z)}{e(z)-1} dz +  \int_{ \om_{ \ep,k}} \frac{ f(z)}{e(z)-1} dz +  \int_{ \ga_{ \ep,k}} \frac{ f(z)}{e(z)-1} dz + \int_{ \la_{ \ep,k}^-} \frac{ f(z)}{e(z)-1} dz $

\end{lemme}

\noindent {\bf Preuve}: \\
On  définit  $ g $ par $ \ds g(z) = \frac{f(z)}{e(z)-1} $, $ g $ est alors méromorphe sur $ U $. \\
Pour $ \ds m \in [[2,k]] \ \ Res(g,m) = \frac{f(m)}{e'(m)} = \frac{f(m)}{2 i \pi} $. Le résultat découle donc du théorème des résidus. \\

\begin{prop}\label{formule de représentation intégrale en dimension N}
Soient $ \ep > 0 $,  $ k \in \N^* $ et  $ U $  un ouvert simplement connexe de $ \C $ contenant  $ R_{\ep , k} $ . \\
Soit $ f:U^N \rightarrow \C $ holomorphe. \\
Pour $ \tau \in {\cal S}_N $ on définit $ f_{\tau} : U^N \to \C $ par $ \ds  f_{\tau}(z_1,...,z_N)=f(z_{\tau(1)},...,z_{\tau(N)}) $ \\
Sous ces conditions $ \ds \sum_{\mb \in [[2,k]]^N } f( \mb ) $ est une somme de $ 4^N $ termes de la forme: 
$$ \int_{( \ga_{\ep,1}^-)^{N_1} \times  (\la_{ \ep,k}^-)^{N_2} \times  (\om_{ \ep,k})^{N_3} \times ( \ga_{ \ep,k})^{N_4}} f_\tau(z_1,...,z_N) \prod _{n=1}^N \frac{1}{e(z_n)-1} d \zb $$
o\`u $ N_1,N_2,N_3,N_4 \in \N $ vérifient $ N_1+N_2+N_3+N_4=N $ et $ \tau \in {\cal S}_N $. 

\end{prop}

\noindent {\bf Preuve}: \\
elle se fait par récurrence sur $ N $ en itérant le lemme \ref{formule de représentation intégrale en dimension 1}. \\

\part*{Preuve du théorème \ref{existence du prolongement holomorphe de Z}}

\begin{notation}
pour $ a,b \in \N $ tels que $ a \leq b $ on note: $ [[a,b]] = \{ m \in \N \ | \ a \leq m \leq b \} $. \\
Pour $ a \in \N $ on note $ [[a, + \infty [[ = \{ m \in \N \ | \ a \leq m \} $. \\ 
\end{notation}

  Le lemme suivant s'inspire fortement de l'analogue se trouvant dans \cite{essoua}.

\begin{lemme}\label{controle près de l'axe}

Soient $ P \in \R[X_1,...,X_N] $ et $ \ep_0 > 0  $ tels que: \\
i) $ \forall \xb \in J^N \ P( \xb) >0 $ \\
ii) $ \forall \ab \in \N^N, \al_n \geq 1 \Rightarrow \p^{\ab} P( \xb) \ll x_n ^{-\ep_0}P(\xb) \ \ ( \xb \in J^N ) $ \\
Alors: \\
$ \exists \ \ep>0 $ tel que: \\
i') $ \xb \in J^N $ et $ \yb  \in [-2 \ep, 2 \ep] ^N \Rightarrow  \Re(P( \xb + i \yb )) \geq \frac{1}{2} P(\xb) $ \\
ii') $ \forall \ab \in \N^N, \al_n \geq 1 \Rightarrow \p^{\ab} P( \xb  + i \yb ) \ll x_n ^{-\ep_0}P(\xb) \  \ ( \xb \in J^N  \  \yb  \in [-2 \ep, 2 \ep] ^N ) $. 

\end{lemme}

\noindent {\bf Preuve}: \\
on note $ p = deg(P) $. \\
La formule de Taylor s'écrit: 
$ \ds P(\xb+ i \yb) = \sum_{|\ab| \leq p } \frac{(i \yb) ^{\ab}}{ \ab !} \p^{\ab} P(\xb) = P( \xb) +  \sum_{ 0 < |\ab| \leq p } \frac{(i \yb) ^{\ab}}{ \ab !} \p^{\ab} P(\xb) $ \\
De l'hypothèse ii) on déduit qu'il existe $ c> 0 $ tel que: $ \ds \forall \ab \in \N^N \setminus \{\zerob\} \ \ \forall \xb \in J^N \ |\p^{\ab} P(\xb)| \leq c P(\xb) $. \\
On pose $ \ds A = c \sum_{ 0 < |\ab| \leq p } \frac{1}{ \ab !} $ . \\
Fixons $ \ds 0< \ep \leq \frac{1}{2} $.\\
On a: $ \yb \in [-2 \ep, 2 \ep] ^N \ \ab \neq \zerob  \Rightarrow |\yb^{\ab}| \leq 2 \ep $. \\
De ce qui précède on déduit: $ \forall \xb \in J^N \ \forall \yb  \in [-2 \ep, 2 \ep] ^N \ \ | P(\xb+ i \yb) -  P(\xb) | \leq 2 \ep A P( \xb) $ \\
On pose $ \ds \ep = \frac{1}{4A+2} $; on a alors: $ \ds \forall \xb \in J^N \ \forall \yb  \in [-2 \ep, 2 \ep] ^N \ | P(\xb+ i \yb) -  P(\xb) | \leq \frac{1}{2} P( \xb) $. \\
Il vient: \\
$ \ds \Re( P(\xb+ i \yb)) = P(\xb) +  \Re( P(\xb+ i\yb) - P(\xb)) \geq  P(\xb) - | P(\xb+ i\yb)) - P(\xb)|  \geq  \frac{1}{2}P(\xb) $ \\
Soit $ \ab \in \N^N $ tel que $ \al_n \geq 1 $. Alors:
$$ (\p^{\ab} P)(\xb+ i \yb)) =   \sum_{ |\bb| \leq p } \frac{(i \yb) ^{\bb}}{ \bb !} \p^{\ab + \bb} P(\xb)  \ll x_n ^{-\ep_0}P(\xb) \ \ ( \xb \in J^N  \ \yb  \in [-2 \ep, 2 \ep] ^N ) $$

\begin{lemme}\label{découpage de N^N}
Soit $ N \in \N^* $. \\
Alors $ \N^{*N} $ peut se partitionner ainsi: \\
$ \ds \N^{*N} = \bigsqcup_{j=1}^J A_j $ o\`u $ \ \forall j \ A_j $ est de la forme $ \ds \prod_{n=1}^N B_n $ avec $ B_n = \{1\} $ ou $ B_n = [[2, + \infty [[ $. \\
\end{lemme}

\noindent {\bf Preuve}: \\
elle se fait par récurrence sur $ N \geq 1 $. \\

\noindent{\bf Preuve du théorème \ref{existence du prolongement holomorphe de Z}}

\noindent La preuve se décompose en 2 étapes.

\noindent Etape 1: $ \ds \sb \mapsto Z^*(\sb) \overset{ \text{déf}}{=} \sum_{\mb \in [[2,+\infty[[^N} \mub^{\mb} Q(\mb) \prod_{t=1}^T P_t(\mb)^{- s_t} $ possède un prolongement holomorphe à $ \C^T $. \\
\noindent Preuve de l'étape 1. \\

Les hypothèses permettent de choisir $ \ep_0> 0 $ tel que: 

$ \ds \star \ \prod_{t=1}^T P_t(\xb) \gg \left( \prod_{n=1}^N x_n \right)^{\ep_0} \ \ ( \xb \in J^N ) $ 

$ \ds  \star \ \ab \in \N^N  \ \al_n \geq 1 \Rightarrow \frac{\p^{\ab} P_t}{P_t}(\xb) \ll x_n^{-\ep_0} \ \  ( \xb \in J^N ) $.

\noindent Pour $ 1 \leq n \leq N $, on écrit $ \mu_n = e^{i \th_n} $ o\`u $ \th_n \in \R \setminus 2 \pi \Z $.
Avec cette notation on a:
$$ \sum_{\mb \in [[2,+\infty[[^N} \mub^{\mb} Q(\mb) \prod_{t=1}^T P_t(\mb)^{- s_t}  = \sum_{\mb \in \N^{*N}} Q(\mb) \prod_{n=1}^N e^{i \th_n m_n} \prod_{t=1}^T P_t(\mb)^{- s_t} $$
Pour $ 1 \leq t \leq T $ on applique le lemme \ref{controle près de l'axe} à $ P_t $, ce qui fournit $ \ep_t >0 $. \\
On pose $ \ep = \min \{ \ep_t \ | \ 1 \leq t  \leq T \} $, $ \ep $ est alors fixé pour toute la preuve. \\
Pour $ \sb \in \C^N $, on définit $ f_{\sb}: (]1, + \infty[ + i]-2 \ep, 2 \ep [)^N \to \C $ par $ \ds f_{\sb}(\zb) = Q(\zb) \prod_{t=1}^T P_t(\zb) ^{-s_t} \prod_{n=1}^N e^{i \th_n z_n} $. \\
Gr\^ace au choix de $ \ep $ ceci a un sens et $ f_{\sb} $ est holomorphe. \\ 
Pour $ k \in \N^* $ on a: $ \ds  \sum_{\mb \in [[2,k]]^N} Q(\mb) \prod_{n=1}^N e^{i \th_n m_n} \prod_{t=1}^T P_t(\mb)^{- s_t} = \sum_{\mb \in [[2,k]]^N } f_{\sb}(\mb) $ \\

$ R_{\ep,k} \subset ]1, + \infty[ + i]-2 \ep, 2 \ep [ $ donc on peut appliquer la proposition \ref{formule de représentation intégrale en dimension N} à $ f_{\sb} $, ce qui permet d'écrire $ \ds \sum_{\mb \in [[2,k]]^N } f_{\sb}(\mb) $ comme une somme de $ 4^N $ intégrales. On va s'occuper de celles pour lesquelles $ \ds \tau = Id_{[[1,N]]} $, les autres se traiteraient exactement de la m\^eme manière. \\
On se limite donc à des expressions de la forme: 
$$ \int_{( \ga_{\ep,1}^-)^{N_1} \times  (\la_{ \ep,k}^-)^{N_2} \times  (\om_{ \ep,k})^{N_3} \times ( \ga_{ \ep,k})^{N_4}} Q(\zb) \prod_{t=1}^T P_t(\zb) ^{-s_t} \prod _{n=1}^N \frac{\exp(i \th_n z_n)}{e(z_n)-1} d \zb $$
soit encore:
$$ (-1) ^{N_1 + N_2} \int_{( \ga_{\ep,1})^{N_1} \times  (\la_{ \ep,k})^{N_2} \times  (\om_{ \ep,k})^{N_3} \times ( \ga_{ \ep,k})^{N_4}} Q(\zb) \prod_{t=1}^T P_t(\zb) ^{-s_t} \prod _{n=1}^N \frac{\exp(i \th_n z_n)}{e(z_n)-1} d \zb $$
o\`u $ N_1,N_2,N_3,N_4 \in \N $ vérifient $ N_1+N_2+N_3+N_4=N $. \\
Si $ N_4 \geq 1 $ et si $ \si_1,...,\si_T $ sont assez grands, on montre par convergence dominée que cette expression tend vers $ 0 $ quand $ k $ tend vers $ + \infty $. \\
Si $ N_4 = 0 $ et si $ \si_1,...,\si_T $ sont assez grands, on montre par convergence dominée que, lorsque $ k $ tend vers $ + \infty $, cette expression tend vers:
$$ Y^{N_1,N_2,N_3}(\sb) \overset{ \text{déf}}{=}  (-1) ^{N_1 + N_2} \int_{( \ga_{\ep,1})^{N_1} \times  (\la_{ \ep})^{N_2} \times  (\om_{ \ep})^{N_3} } Q(\zb) \prod_{t=1}^T P_t(\zb) ^{-s_t} \prod _{n=1}^N \frac{\exp(i \th_n z_n)}{e(z_n)-1} d \zb $$ 
On a donc montré qu'il existe $ r > 0 $ tel que sur $ \{ \sb \in \C^T \ | \ \si_1,...,\si_T > r \} $  $ Z^* $ est une combinaison linéaire d'intégrales de la forme $ Y^{N_1,N_2,N_3} $ à permutation près. \\

\noindent Pour conclure il nous suffit donc de montrer que $ Y^{N_1,N_2,N_3} $ possède un prolongement holomorphe à $ \C^T $. \\ \\
$ \star $ Pour $ 1 \leq n \leq N_1 $ on définit $ f_n: [-1,1] \to \C $ ainsi: \\
$ \ds f_n(x) = \frac{\exp(i \th_n \ga_{\ep,1}(x))}{e(\ga_{\ep,1}(x)) - 1} = \frac{\exp(i \th_n(\frac{3}{2} + i \ep x))}{\exp(2i\pi (\frac{3}{2} + i \ep x)) - 1} = - \exp \left( \frac{3}{2} i \th_n \right) \frac{\exp(- \ep \th_n x)}{\exp(- 2\pi \ep x) + 1}  $  \\
$ f_n $ est manifestement continue. \\
$ \star $ Pour $ N_1 + 1 \leq n \leq N_1 + N_2 $ on définit $ \ds f_n : \left[\frac{3}{2}, +\infty \right[ \to \C $ ainsi: \\
$ \ds f_n(x) = \frac{\exp(i \th_n \la_{\ep}(x))}{e(\la_{\ep}(x)) - 1} = \frac{\exp(i \th_n(x + i \ep))}{\exp(2i\pi (x + i \ep )) - 1} = - \exp(- \ep \th_n) \frac{\exp(i \th_n x)}{1 - \exp(- 2 \pi \ep )\exp( i 2\pi x ) }  $  \\
$ \ds \frac{\th_n}{2 \pi} \notin \Z $ donc (voir l'exemple \ref{exemples de fonctions appartenant à B(r)}) $ \ds f_n \in \br \left( \frac{3}{2} \right) $. \\
$ \star $ Pour $ N_1 + N_2 + 1 \leq n \leq N $ on définit $ \ds f_n : \left[\frac{3}{2}, +\infty \right[ \to \C $ ainsi: \\
$ \ds f_n(x) = \frac{\exp(i \th_n \om_{\ep}(x))}{e(\om_{\ep}(x)) - 1} = \frac{\exp(i \th_n (x - i \ep))}{\exp(2i \pi (x - i \ep )) - 1} = - \exp( \ep \th_n) \frac{\exp(i \th_n x)}{1 - \exp(2 \pi \ep ) \exp( i 2 \pi x ) }  $  \\
$ \ds \frac{\th_n}{2 \pi} \notin \Z $ donc (voir l'exemple \ref{exemples de fonctions appartenant à B(r)}) $ \ds f_n \in \br \left(\frac{3}{2} \right) $. \\ \\
Pour $ P \in \C[X_1,...,X_N] $ et $ N_1,N_2,N_3 $ de somme $ N $, on définit $ P^{N_1,N_2,N_3} \in \C[X_1,...,X_N] $ par:
\begin{align*}
P^{N_1,N_2,N_3}(\xb) & = P(\ga_{\ep,1}(x_1),...,\ga_{\ep,1}(x_{N_1}),\la_{\ep}(x_{N_1 + 1}),...,\la_{\ep}(x_{N_1 + N_2}),\om_{\ep}(x_{N_1 + N_2 + 1}),...,\om_{\ep}(x_N)) \\
                     & = P(\frac{3}{2} + i \ep x_1,...,\frac{3}{2} + i\ep x_{N_1},x_{N_1 + 1} + i \ep,...,x_{N_1 + N_2} + i \ep,x_{N_1 + N_2 + 1} - i \ep,...,x_N - i \ep) 
\end{align*}
Muni de ces notations, on constate que:
$$ Y^{N_1,N_2,N_3}(\sb) =  (-1)^{N_1 + N_2} (i \ep)^{N_1} \int_{ [-1,1]^{N_1} \times [\frac{3}{2}, + \infty[^{N - N_1}} Q^{N_1,N_2,N_3}(\xb) \prod_{t=1}^T P_t^{N_1,N_2,N_3}(\xb)^{-s_t} \prod_{n=1}^N f_n(x_n) d \xb $$
Il nous suffit maintenant de vérifier les hypothèses du théorème \ref{existence du prolongement holomorphe de Y} (que l'on applique ici sur \\
$ \ds [-1,1]^{N_1} \times \left[\frac{3}{2}, + \infty \right[^{N - N_1} $ et non sur $ [-1,1]^{N_1} \times J^{N - N_1} $, ce qui ne pose clairement aucun problème). \\
$ \ast \ f:[-1,1]^{N_1} \to \C $ définie par $ \ds f(x_1,...,x_{N_1}) = \prod_{n=1}^{N_1} f_n(x_n) $ est clairement continue. \\
$ \ast $ On a déjà vu que  $ \ds f_{N_1 + 1},...,f_N  \in \br \left(\frac{3}{2} \right) $. \\
$ \ds \ast \ P_t^{N_1,N_2,N_3}(\xb) = P_t \left( \left(\frac{3}{2},...,\frac{3}{2},x_{N_1 + 1},...,x_N \right) + i(\ep x_1,..., \ep x_{N_1},\ep,...,\ep,-\ep,...,-\ep) \right) $ \\
or on a déterminé $ \ep $ gr\^ace au lemme \ref{controle près de l'axe}, donc: \\
$ \ds \forall \xb \in [-1,1]^{N_1} \times J^{N - N_1} \ \ \Re \left(P_t^{N_1,N_2,N_3}(x_1,...,x_N) \right) \geq \frac{1}{2} P_t \left(\frac{3}{2},...,\frac{3}{2},x_{N_1 + 1},...,x_N \right) $ \\
On  déduit de ceci que $ \ds \forall \xb \in [-1,1]^{N_1} \times \left[ \frac{3}{2},+\infty \right[^{N-N_1}$ on a:

$ \ds \star \ \Re \left( P_t^{N_1,N_2,N_3}(\xb) \right) > 0 $ 

$ \ds \star \ \left| P_t^{N_1,N_2,N_3}(\xb) \right| \geq \frac{1}{2} P_t \left( \frac{3}{2},...,\frac{3}{2},x_{N_1 + 1},...,x_N \right) $ \\
De cette dernière inégalité, on déduit toute suite: \\
$ \ds \left| P_t^{N_1,N_2,N_3}(\xb) \right| \gg 1 \ \ \left( \xb \in [-1,1]^{N_1} \times \left[ \frac{3}{2},+\infty \right[^{N-N_1} \right) $ \\
et $ \ds \forall \xb \in [-1,1]^{N_1} \times \left[ \frac{3}{2},+\infty \right[^{N-N_1} \ \ \ \prod_{t=1}^T \left| P_t^{N_1,N_2,N_3}(\xb) \right| \geq \left( \frac{1}{2} \right)^T \prod_{t=1}^T P_t \left( \frac{3}{2},...,\frac{3}{2},x_{N_1 + 1},...,x_N \right) $ \\
Il vient: \\ $ \ds \prod_{t=1}^T \left| P_t^{N_1,N_2,N_3}(\xb) \right| \gg \left( \prod_{n= N_1 + 1}^N x_n \right)^{\ep_0} \ \ \left( \xb \in [-1,1]^{N_1} \times \left[ \frac{3}{2},+ \infty \right[^{N - N_1} \right) $ \\
Si $ \ab \in { \{  0 \} }^{N_1} \times \N^{N - N_1} $ et $ N_1 + 1 \leq n \leq N  $ sont tels que $ \al_n \geq 1 $, alors:
\begin{align*}
\p^{\ab} P_t^{N_1,N_2,N_3}(\xb) & =  (\p^{\ab} P_t) \left( \frac{3}{2} + i \ep x_1,...,\frac{3}{2} + i\ep x_{N_1},x_{N_1 + 1} + i \ep,...,x_{N_1 + N_2} + i \ep,x_{N_1 + N_2 + 1} - i \ep,...,x_N - i \ep \right) \\
         & \ \ll x_n^{ - \ep_0}  P_t \left( \frac{3}{2},...,\frac{3}{2},x_{N_1 + 1},...,x_N \right)  \ \  \left( \xb \in [-1,1]^{N_1} \times \left[ \frac{3}{2},+ \infty \right[^{N-N_1} \right) \\
         & \ \ll x_n^{ - \ep_0} \left| P_t^{ N_1,N_2,N_3}(\xb) \right|  \ \  \left( \xb \in [-1,1]^{N_1} \times \left[ \frac{3}{2},+ \infty \right[^{N-N_1} \right)
\end{align*}
Ainsi s'achève les vérifications des hypothèses du théorème \ref{existence du prolongement holomorphe de Y}, ce qui termine la preuve de l'étape 1. \\

\noindent Etape 2: conclusion. \\ \\
On va montrer le théorème \ref{existence du prolongement holomorphe de Z} par récurrence sur $ N \geq 1 $. \\
$ \star $ Pour $ N = 1 $, il suffit d'écrire:
$$ Z(Q,P_1,...,P_T,\mu,\sb) = \mu Q(1) \prod_{t=1}^T P_t(1)^{- s_t} + \sum_{m \geq 2} \mu^m Q(m) \prod_{t=1}^T P_t(m)^{- s_t} $$
L'étape 1 permet alors de conclure. \\
$ \star $ Si le résultat est vrai pour tout $ n $ compris entre $ 1 $ et $ N - 1 $, alors gr\^ace au lemme \ref{découpage de N^N} et à l'étape 1, on voit qu'il est vrai pour $ N $.

\part*{Lemme d'échange et valeurs aux points de $ (-\N)^T $}

\begin{prop}\label{préliminaire au lemme d'échange}
Soient $ Q, P_1,...,P_T \in \R[X_1,...,X_N] $ et $ 1 \leq T_0 \leq T - 1 $. \\
On suppose que: \\
a) $ P_1,...,P_T $ vérifient HR \\
b) $ \ds \prod_{t=1}^{T_0} P_t(\xb) \xrightarrow[ \substack{ |\xb| \rightarrow + \infty \\ \xb \in J^N}]{} + \infty $. \\
Soient de plus $ \mub \in (\T \setminus \{1\})^N $ et $ k_1,...,k_T \in \N $. \\
Alors:
$$ Z(Q,P_1,...,P_T,\mub, -k_1,...,-k_T) = Z \left(Q \prod_{t = T_0 +1}^T {P_t}^{k_t},P_1,...,P_{T_0},\mub,-k_1,...,-k_{T_0} \right) $$
\end{prop}

\noindent {\bf Preuve}: \\
On définit $ f: \C^{T_0} \to \C $ par $ f(s_1,...,s_{T_0}) =  Z(Q,P_1,...,P_T,\mub,s_1,...,s_{T_0},-k_{T_0 + 1},...,-k_T) $. \\
$ f $ est holomorphe. \\
Par la proposition \ref{pour sigma grand on est dans le domaine de convergence de Z} il existe $ \si_0 \in \R $ tel que pour $ \si_1,...,\si_{T_0} \geq \si_0 $ on ait:
$$ Z(Q,P_1,...,P_T,\mub,\si_1,...,\si_{T_0},-k_{T_0 + 1},...,-k_T) = \sum_{\mb \in \N^{*N}} \mub^{\mb} Q(\mb) \prod_{t=1}^{T_0} P_t(\mb)^{- \si_t} \prod_{t = T_0 + 1}^T P_t(\mb)^{k_t}  $$ 
On définit $ g: \C^{T_0} \to \C $ par $ \ds g(s_1,...,s_{T_0}) = Z \left( Q \prod_{t = T_0 +1}^T {P_t}^{k_t},P_1,...,P_{T_0},\mub,s_1,...,s_{T_0} \right) $. \\
$ g $ est holomorphe. \\
Par la proposition \ref{pour sigma grand on est dans le domaine de convergence de Z} il existe $ \si_0' \in \R $ tel que pour $ \si_1,...,\si_{T_0} \geq \si_0' $ on ait:
$$ Z \left( Q \prod_{t = T_0 +1}^T {P_t}^{k_t},P_1,...,P_{T_0},\mub,s_1,...,s_{T_0} \right) = \sum_{\mb \in \N^{*N}} \mub^{\mb} Q(\mb) \prod_{t = T_0 +1}^T {P_t}^{k_t} \prod_{t=1}^{T_0} P_t(\mb)^{- \si_t}  $$
On constate que pour  $ \si_1,...,\si_{T_0} \geq \max(\si_0,\si_0') $ on a $ f(\si_1,...,\si_{T_0}) =  g(\si_1,...,\si_{T_0}) $; par prolongement analytique on en déduit que $ f = g $. \\ 
En particulier $ f(-k_1,...,-k_{T_0}) = g(-k_1,...,-k_{T_0}) $, ce qui est exactement le résultat voulu. \\

\noindent {\bf Preuve du lemme d'échange}: \\
La proposition \ref{préliminaire au lemme d'échange} permet d'affirmer que les quantités considérées sont toutes deux égales à: \\
$ Z(Q,P_1,...,P_T,Q_1,...,Q_{T'},\mub,-k_1,...,-k_T,-l_1,...,-l_{T'}) $. \\

\begin{lemme}\label{valeurs en 0}
Soient $ Q \in \R[X_1,...,X_N] $ et $ \mub \in (\T \setminus \{1\})^N $. \\
On note $ \ds Q = \sum_{\ab \in S} a_{\ab} \Xb^{\ab} $. \\
Alors:
 $$  Z(Q,X_1,...,X_N,\mub,0,...,0) = \sum_{\ab \in S} a_{\ab} \prod_{n=1}^N  \zeta_{\mu_n}(- \al_n)  $$
\end{lemme}

\noindent {\bf Preuve}: \\
si $ \si_1,...,\si_N $ sont suffisamment grands on a: 
\begin{align*}
Z(Q,X_1,...,X_N,\mub,\sb) & = Z \left( \sum_{\ab \in S} a_{\ab} \Xb^{\ab},X_1,...,X_N,\mub,\sb \right) \\
                          & = \sum_{\ab \in S} a_{\ab} Z(\Xb^{\ab},X_1,...,X_N,\mub,\sb) \\
                          & = \sum_{\ab \in S} a_{\ab} \sum_{\mb \in \N^{*N}} \mub^{\mb} \mb^{\ab} \prod_{n=1}^N m_n^{-s_n} \\ 
                          & = \sum_{\ab \in S} a_{\ab} \sum_{m_1,...,m_N \geq 1} \prod_{n=1}^N \mu_n^{m_n} m_n^{\al_n - s_n} \\
                          & = \sum_{\ab \in S} a_{\ab} \prod_{n=1}^N \sum_{m_n \geq 1} \mu_n^{m_n} m_n^{\al_n - s_n} \\ 
                          & = \sum_{\ab \in S} a_{\ab} \prod_{n=1}^N \zeta_{\mu_n}(s_n - \al_n) 
\end{align*}
Par prolongement analytique, on a donc:
$$ \forall \sb \in \C^N \ \ Z(Q,X_1,...,X_N,\mub,\sb) = \sum_{\ab \in S} a_{\ab} \prod_{n=1}^N \zeta_{\mu_n}(s_n - \al_n)  $$
Il suffit maintenant de faire $ \sb = \zerob $ dans cette égalité pour obtenir le résultat cherché. \\

\noindent {\bf Preuve du théorème \ref{formule pour les valeurs de Z aux entiers négatifs} }
\begin{align*}
Z(Q,P_1,...,P_T,\mub,- k_1,...,- k_T) & = Z \left(Q \prod_{n=1}^N X_n^0,P_1,...,P_T,\mub,- k_1,...,- k_T \right) \\
\text{gr\^ace au lemme \ref{lemme d'échange}}   \ \   & = Z \left(Q \prod_{t=1}^T P_t^{k_t},X_1,...,X_N,\mub,0,...,0 \right) \\
\text{gr\^ace au lemme \ref{valeurs en 0}}   \  \  & = \sum_{\ab \in S} a_{\ab} \prod_{n=1}^N \zeta_{\mu_n}(- \al_n) 
\end{align*}

\part*{Une formule pour les valeurs de $ \zeta_{\mu} $ aux entiers négatifs}

Le lemme suivant se trouve dans \cite{zagier}. \\

\begin{lemme}\label{valeurs aux entiers négatifs et dl}
Soit $ (a_m)_{m \in \N^*} $ une suite de nombres complexes. \\
On pose $ \ds Z(s) = \sum_{m=1}^{+\infty} \frac{a_m}{m^s} $ et l'on suppose qu'il existe $ s \in \C $ tel que cette série converge. \\
Grâce à cette hypothèse, on peut définir $ f \colon \R_+^* \to \C $ par: $ \ds f(x) = \sum_{m=1}^{+\infty} a_m e^{-mx} $. \\
On suppose qu'il existe une suite $ (c_k)_{k \in \N} $ de complexes telle que, pour tout $ K \in \N^* $, on ait au voisinage de $ 0 $: $ \ds f(x) = \sum_ {k=0}^{K-1} c_k x^k + O(x^K) $. \\
Alors: $ Z $ se prolonge holomorphiquement à $ \C $ et $ \forall k \in \N \ Z(-k) = (-1)^k k! c_k $. \\
\end{lemme}

Nous aurons besoin des nombres de Stirling de seconde espèce. Rappelons en tout d'abord la définition: \\

\begin{déf}
Soient $ k,\ell \in \N $. \\
Par définition le nombre de Stirling de second espèce (associé à $ (k,\ell) $) est le nombre de partitions en $ \ell $ parties d'un ensemble à $ k $ éléments. Cet entier naturel est noté $ S(k,\ell) $. \\
\end{déf}

\begin{ex}
$ S(0,0) = 1 $; pour $ k \in \N $ $ S(k,k) = 1 $; si $ \ 0 \leq k < \ell $ alors $ S(k,\ell) = 0 $. \\
\end{ex} 

On va maintenant rappeler quelques propriétés élémentaires de ces nombres. Pour les preuves on renvoie par exemple à \cite{comtet}. \\

\begin{lemme}
$ \forall k \in \N \ \forall \ell \in \N^* \ S(k+1,\ell) = \ell S(k,\ell) + S(k,\ell-1) $. \\
\end{lemme}

\begin{lemme}\label{formule explicite pour les nombres de stirling}
Pour tous $ k,\ell \in \N $ on a: $ \ds S(k,\ell) = \frac{1}{\ell !} \sum_{j=0}^\ell (-1)^{\ell - j} \binom{\ell}{j} j^k $. \\
\end{lemme}

\begin{lemme}\label{dérivées de g circ exp}
Soit $ g \colon \R_+^* \to \C $ de classe $ C^{\infty} $. On définit $ f \colon \R \to \C $ par $ f = g \circ \exp $.  \\
Alors pour tout $ k \in \N $ on a: $ \ds \forall x \in \R \ \ f^{(k)}(x) = \sum_{\ell=0}^k S(k,\ell) e^{\ell x} g^{(\ell)}(e^x) $. \\
\end{lemme}

\noindent{\bf Preuve}: \\
elle se fait par récurrence sur $ k \in \N $. \\
$ \star $ Pour $ k=0 $ cela découle de $ S(0,0) = 1 $. \\
$ \star $ Si l'assertion est vraie au rang $ k $, alors pour tout $ x \in \R $ on a:
\begin{align*}
f^{(k+1)}(x) & = \sum_{\ell=0}^k S(k,\ell) \left(\ell e^{\ell x} g^{(\ell)}(e^x) + e^{\ell x} e^x g^{(\ell+1)}(e^x) \right) \\
             & =  \sum_{\ell=0}^k S(k,\ell) \ell e^{\ell x} g^{(\ell )}(e^x) +  \sum_{\ell =1}^{k+1} S(k,\ell -1) e^{\ell x} g^{(\ell )}(e^x) 
\end{align*}
$ S(k,k+1) = 0 $ donc:
\begin{align*}
f^{(k+1)}(x) & = \sum_{\ell =1}^{k+1} \left[ \ell  S(k,\ell ) + S(k,\ell -1) \right] e^{\ell x} g^{(\ell )}(e^x) \\
             & = \sum_{\ell =1}^{k+1} S(k+1,\ell ) e^{\ell x} g^{(\ell )}(e^x) \\
             & = \sum_{\ell =0}^{k+1} S(k+1,\ell ) e^{\ell x} g^{(\ell )}(e^x)
\end{align*}
ce qui termine la récurrence. \\

Nous pouvons maintenant prouver le: \\

\begin{lemme}
Soit $ \mu \in \T \setminus \{1\} $. \\
Alors pour tout $ k \in \N $ on a: $ \ds \zeta_{\mu}(-k) = \frac{(-1)^k \mu}{1-\mu} \sum_{\ell =0}^k \frac{\ell  ! S(k,\ell )}{(\mu - 1)^\ell } $. \\
\end{lemme}

\noindent {\bf Preuve}: 
$$ \forall x>0 \ \sum_{m=1}^{+\infty} \mu^m e^{-mx} =  \sum_{m=1}^{+\infty} \left( \mu e^{-x} \right)^m = \mu e^{-x} \frac{1}{1 - \mu e^{-x}}  = \frac{\mu}{e^x - \mu}   $$
On définit $ f \colon \R \to \C $ par $ \ds f(x) = \frac{\mu}{e^x - \mu} $  et $ g \colon \R_+^* \to \C $ par $ \ds g(y) = \frac{\mu}{y - \mu} $. \\
$ g $ est $ C^{\infty } $ et $ f = g \circ \exp $ donc \ref{dérivées de g circ exp} s'applique et donne:
$$ \forall x \in \R \ \ f^{(k)}(x) = \sum_{\ell =0}^k S(k,\ell ) e^{\ell x} g^{(\ell )}(e^x). $$
$ \ds \forall y \in \R_+^* \ g(y) = -\mu \frac{1}{\mu - y } $ donc: $ \ds \forall \ell  \in \N\ \forall y \in \R_+^* \ g^{(\ell )}(y) = - \mu \frac{\ell !}{(\mu - y)^{\ell+1} } $. \\ 
De ce qui précède on déduit que: $ \ds f^{(k)}(0) = \sum_{\ell=0}^k S(k,\ell) \left(- \mu \frac{\ell!}{(\mu - 1)^{\ell+1} } \right) $, \\
on peut maintenant conclure en utilisant \ref{valeurs aux entiers négatifs et dl}. \\

\part*{Preuve des théorèmes 4 et 5}

\noindent {\bf Preuve}: \\
soit $ \kb \in \N^T $. On note $ S_{\kb} $ le support de $ \ds Q \prod_{t=1}^T P_t^{k_t} $. Soit $ (a_{\ab})_{\ab \in S_{\kb}} $ telle que $ \ds \sum_{\ab \in S_{\kb}} a_{\ab} \Xb^{\ab} = Q \prod_{t=1}^T P_t^{k_t} $. \\
Le théorème B dit que $ \ds Z(Q,P_1,...,P_T, \mub,-\kb) = \sum_{\ab \in S_{\kb}} a_{\ab} \prod_{n=1}^N \zeta_{\mu_n}(-\al_n) $. \\
On va maintenant utiliser le lemme de la section suivante  sous la forme suivante:
$$ \forall k \in \N \ \ \zeta_{\mu}(-k) = \frac{(-1)^k \mu}{1-\mu} \sum_{\ell \in \N} \frac{\ell ! S(k,\ell)}{(\mu - 1)^\ell}. $$
Dans cette formule la somme est en réalité une somme finie. \\
Dans le calcul qui suit toutes les sommes sont en réalité des sommes finies. 
\begin{align*}
Z(Q,P_1,...,P_T, \mub,-\kb) & = \sum_{\ab \in S_{\kb}} \left[ a_{\ab} \prod_{n=1}^N \left( \frac{(-1)^{\al_n} \mu_n}{1-\mu_n} \sum_{\ell_n \in \N} \frac{\ell_n ! S(\al_n,\ell_n)}{(\mu_n - 1)^{\ell_n}} \right) \right] \\ 
                           & = \frac{ \mub^{\unb}}{(\unb - \mub)^{\unb}} \sum_{\ab \in S_{\kb}}\left[ (-1)^{|\ab|} a_{\ab} \prod_{n=1}^N \sum_{\ell_n \in \N} \frac{\ell_n ! S(\al_n,\ell_n)}{(\mu_n - 1)^{\ell_n}} \right] \\
                           & = \frac{ \mub^{\unb}}{(\unb - \mub)^{\unb}} \sum_{\ab \in S_{\kb}}\left[ (-1)^{|\ab|} a_{\ab} \sum_{\lb \in \N^N} \prod_{n=1}^N \frac{\ell_n ! S(\al_n,\ell_n)}{(\mu_n - 1)^{\ell_n}} \right] \\
                           & = \frac{ \mub^{\unb}}{(\unb - \mub)^{\unb}} \sum_{\lb \in \N^N} \left\{ \prod_{n=1}^N \left( \frac{\ell_n !}{(\mu_n - 1)^{\ell_n}} \right) \sum_{\ab \in S_{\kb}} \left( (-1)^{|\ab|} a_{\ab} \prod_{n=1}^N S(\al_n,l_n) \right) \right\}
\end{align*}
Pour $ \lb \in \N^N $, on définit $ Z_{\lb} \colon (-\N)^T \to \Q $ par:
$$ \forall \kb \in \N^T \  Z_{\lb}(-\kb) = \prod_{n=1}^N \left( \frac{\ell_n !}{(\mu_n - 1)^{\ell_n}} \right)  \sum_{\ab \in S_{\kb}} \left( (-1)^{|\ab|} a_{\ab} \prod_{n=1}^N S(\al_n,\ell_n)  \right) $$
et l'on remarque que $ \ds \forall \kb \in \N^T \ |Z_{\lb}(-\kb)|_p \leq \prod_{n=1}^N \frac{|\ell_n !|_p}{|\mu_n - 1|_p^{\ell_n}} $. \\
En utilisant \ref{formule explicite pour les nombres de stirling} il vient: 
\begin{align*}
Z_{\lb}(-\kb) & = \prod_{n=1}^N \left( \frac{\ell_n !}{(\mu_n - 1)^{\ell_n}} \right) \sum_{\ab \in S_{\kb}} \left\{ (-1)^{|\ab|} a_{\ab} \prod_{n=1}^N \left[ \frac{1}{\ell_n !} \sum_{j_n=0}^{\ell_n} \left( (-1)^{\ell_n-j_n} \binom{\ell_n}{j_n} j_n^{\al_n} \right) \right] \right\}\\
              & = (\unb - \mub)^{-\lb} \sum_{\ab \in S_{\kb}} \left\{ (-1)^{|\ab|} a_{\ab} \prod_{n=1}^N \left[ \sum_{j_n=0}^{\ell_n} \left( (-1)^{j_n} \binom{\ell_n}{j_n} j_n^{\al_n} \right) \right] \right\} \\
              & = (\unb - \mub)^{-\lb} \sum_{\ab \in S_{\kb}} \left\{ (-1)^{|\ab|} a_{\ab} \sum_{ \jb \in \prod_{n=1}^N \{0,...,\ell_n \} } \left[ \prod_{n=1}^N \left( (-1)^{j_n} \binom{\ell_n}{j_n} j_n^{\al_n} \right) \right] \right\} \\
              & = (\unb - \mub)^{-\lb} \sum_{\jb \in \prod_{n=1}^N \{0,...,\ell_n \}} \left\{ \prod_{n=1}^N \left[ (-1)^{j_n} \binom{\ell_n}{j_n} \right] \sum_{\ab \in S_{\kb}} \left[ (-1)^{|\ab|} a_{\ab} \prod_{n=1}^N j_n^{\al_n} \right] \right\} \\
              & = (\unb - \mub)^{-\lb} \sum_{\jb \in \prod_{n=1}^N \{0,...,\ell_n \}} \left\{ \prod_{n=1}^N \left[ (-1)^{j_n} \binom{\ell_n}{j_n} \right]  \sum_{\ab \in S_{\kb}} \left[  a_{\ab} \prod_{n=1}^N (-j_n)^{\al_n} \right] \right\} \\
              & = (\unb - \mub)^{-\lb} \sum_{\jb \in \prod_{n=1}^N \{0,...,\ell_n \}} \left\{ (-1)^{|\jb|} \binom{\lb}{\jb} Q(-\jb) \prod_{t=1}^T P_t(-\jb)^{k_t} \right\}
\end{align*}
Pour $ x \in \Z_p^* $ on note $ w(x) $ le teichmüller de $ x $ et $ \ds <x> = \frac{x}{w(x)} $. \\
Soient $ t \in \{1,...,T \} $ et $ \jb \in \N^N $. \\
On remarque que si $ k_t \in \N $ vérifie $ k_t = r_t \ [p-1] $ alors $ \ds P_t(-\jb)^{k_t} = w \left(P_t(-\jb) \right)^{r_t} <P_t(-\jb)>^{k_t} $. \\
Cela nous incite à définir $ Z_{\lb}^{\rb} \colon \Z_p^T \to \C_p $ par: 
$$ Z_{\lb}^{\rb}(s_1,...,s_T) = (\unb - \mub)^{-\lb} \sum_{\jb \in \prod_{n=1}^N \{0,...,\ell_n \}} (-1)^{\jb} \binom{\lb}{\jb} Q(-\jb) \prod_{t=1}^T  w \left(P_t(-\jb) \right)^{r_t} <P_t(-\jb)>^{-s_t} $$
Soit $ \kb \in \N^T $ vérifiant $ \forall t \in \{1,...,T \} \ k_t = r_t \ [p-1] $. Alors: \\
$ Z_{\lb}^{\rb}(-\kb) = Z_{\lb}(-\kb) $ et donc, grâce à une remarque précédente: $ \ds |Z_{\lb}^{\rb}(-\kb)|_p \leq \prod_{n=1}^N \frac{|\ell_n !|_p}{|\mu_n - 1|_p^{\ell_n}} $. \\Comme $ Z_{\lb}^{\rb} $ est continue et que $ \ds - \prod_{t=1}^T \left( r_t + (p-1) \N \right) $ est dense dans $ \Z_p^T $ on en déduit que: 
$$ \forall \sb \in \Z_p^T \ \ |Z_{\lb}^{\rb}(\sb)|_p \leq \prod_{n=1}^N \frac{|\ell_n !|_p}{|\mu_n - 1|_p^{\ell_n}} $$
$ \ds  \frac{|\ell!|_p}{|\mu_n - 1|_p^\ell} \xrightarrow[ \ell \rightarrow + \infty]{} 0 $ donc la définition suivante a un sens: \\
on définit $ Z_p^{\rb}(Q,P_1,...,P_T,\mub,\cdot) \colon \Z_p^T \to \C_p $ par: $ \ds Z_p^{\rb}(Q,P_1,...,P_T,\mub,\sb) = \frac{ \mub^{\unb}}{(\unb - \mub)^{\unb}} \sum_{\lb \in \N^T} Z_{\lb}^{\rb}(\sb) $. \\
Ceci convient clairement. \\

\noindent {\bf Remerciements}: \\
Durant l'élaboration de ce travail j'ai été encadré par Driss Essouabri. Son encadrement fut de grande qualité, tant sur le plan scientifique qu'humain. Qu'il en soit ici remercié! \\
Je tiens aussi à remercier Ben Lichtin pour sa lecture attentive et constructive de ce travail.

\end{document}